\documentclass[a4paper,10pt]{amsart}
\usepackage[margin=1.25in]{geometry}
\usepackage{amssymb,bm,amsfonts,amsmath, psfrag,multicol,amsthm}
\usepackage{tensor}
\usepackage[dvipsnames]{xcolor}
\usepackage[nosort]{cite}
\usepackage{cite}
\usepackage{ytableau}
\usepackage{tikz}
\usetikzlibrary{decorations.pathreplacing,calligraphy,positioning, shapes.geometric}
\usetikzlibrary{arrows,patterns.meta}
\usetikzlibrary{decorations.pathreplacing}
\usepackage{graphicx}
\usepackage{ifpdf}
\usepackage{dashrule}
\usepackage{booktabs}
\usepackage{pdfpages}
\usepackage{epsfig}
\usepackage{epstopdf}
\usepackage{mathrsfs}
\usepackage{stmaryrd}
\usepackage{enumerate}
\usepackage{soul}
\usepackage[colorlinks,urlcolor=RedViolet,linkcolor=red,anchorcolor=green,citecolor=blue]{hyperref}
\usepackage{longtable}
\usepackage{booktabs}
\usepackage{threeparttable}
\usepackage{array}
\usepackage{multirow}
\usepackage{boldline}
\allowdisplaybreaks
\numberwithin{equation}{section}
\newtheorem{thm}{Theorem}[section]

\newtheorem{lem}[thm]{Lemma}

\newtheorem{definition}{Definition}[section]
\newtheorem{prop}{Proposition}[section]

\newtheorem{exam}{Example}[section]

\newtheorem{remark}{Remark}[section]

\def\P{{\mathcal P}}
\def\D{{\mathcal D}}
\def\G{{\mathcal G}}
\def\qed{\hfill \rule{4pt}{7pt}}
\def\pf{\noindent {\it{Proof.} \hskip 2pt}}
\def\ij{\langle i,j \rangle}
\def\babyblue{\textcolor{babyblue}}
\def\red{\textcolor{red}}
\def\des{\operatorname{des}}
\def\ides{\operatorname{ides}}
\def\asc{\operatorname{asc}}
\def\iasc{\operatorname{iasc}}

\definecolor{babyblue}{rgb}{0.54, 0.81, 0.94}
\definecolor{airforceblue}{rgb}{0.36, 0.54, 0.66}
\definecolor{gainsboro}{rgb}{0.86, 0.86, 0.86}
\definecolor{afblue}{rgb}{0.36, 0.54, 0.66}
\definecolor{royalblue}{rgb}{0.25,0.41,0.87}
\def\blue{\textcolor{blue}}
\def\red{\textcolor{red}}
\def\gainsboro{\textcolor{gainsboro}}
\def\afblue{\textcolor{afblue}}
\def\royalblue{\textcolor{royalblue}}

\newcommand{\cpfthm}[1]{\noindent{\emph{Combinatorial Proof of Theorem #1.}\hskip 2pt}}
\newcommand{\cpfcor}[1]{\noindent{\emph{Combinatorial Proof of Corollary #1.}\hskip 2pt}}
\newcommand{\cpflem}[1]{\noindent{\emph{Combinatorial Proof of Lemma #1.}\hskip 2pt}}

\begin{document}

\title[Combinatorial proofs on the descents and inverse descents]{Combinatorial proofs on the joint distribution \\[3pt] of descents and inverse descents}

\author[F.Z.K. Li]{Frank Z.K. Li}
\author[X. Liu]{Xunhao Liu}

\date{\today}



\begin{abstract}
  Let $A_{n,i,j}$ be the number of permutations on $[n]$
with $(i-1)$ descents and $(j-1)$ inverse descents.
Carlitz, Roselle and Scoville in 1966
first revealed some combinatorial and arithmetic properties of $A_{n,i,j}$,
which contain a recurrence of $A_{n,i,j}$.
Using the idea of balls in boxes,
Petersen gave a combinatorial interpretation for the generating function of $A_{n,i,j}$,
and obtained the same recurrence of $A_{n,i,j}$ from its generating function.
Subsequently, Petersen asked whether there is a visual way to understand this recurrence.
In this paper,
after observing the internal structures of permutation grids,
we present a combinatorial proof for the recurrence of $A_{n,i,j}$.
Let $I_{n,k}$ and $J_{n,k}$ be
the number of involutions and fixed-point free involutions on $[n]$
with $k$ descents, respectively.
With the help of algebraic method on generating functions,
Guo and Zeng derived two recurrences of $I_{n,k}$ and $J_{2n,k}$
that play an essential role in the proof of their unimodal properties.
Surprisingly,
the constructive approach to the recurrence of $A_{n,i,j}$
is found to fuel the combinatorial interpretations of
these two recurrences of $I_{n,k}$ and $J_{2n,k}$.
\end{abstract}

\maketitle

\section{Introduction}
Let $n\geq1$ be an integer and $[n]$ be the set $\{1,2,\ldots,n\}$,
then a \emph{permutation} $\pi=\pi_1\pi_2\ldots \pi_n$ on $[n]$ is a bijection from $[n]$ onto itself defined by $\pi(i)=\pi_i$ for $1\leq i\leq n$.
It is known that the set of permutations on $[n]$, denoted by $\mathfrak{S}_n$,
is the \emph{symmetric group}.

For a permutation $\pi=\pi_1\pi_2\ldots \pi_n$,
let $i\in [n-1]$,
then we call $i$ a \emph{descent} if $\pi_i>\pi_{i+1}$
and an \emph{inverse descent} if $\pi^{-1}_i>\pi^{-1}_{i+1}$,
where $\pi^{-1}$ is the inverse of $\pi$ in $\mathfrak{S}_n$.
Correspondingly, we say $i\in [n-1]$ is an \emph{ascent} if $\pi_i<\pi_{i+1}$,
and an \emph{inverse ascent} if $\pi^{-1}_i<\pi^{-1}_{i+1}$.
In the rest of this paper,
we follow the notations in \cite{Rawlings-1984} to abbreviate inverse descents and inverse ascents by \emph{idescents} and \emph{iascents}, respectively.
Given a permutation $\pi$, define the \emph{descent set} $D(\pi)=\{i\mid \pi_i>\pi_{i+1}\}$,
and let $\des(\pi)=|D(\pi)|$.
To unify the notations, we use $\ides(\pi)$ instead of $\des(\pi^{-1})$.
If $\pi=264135$,
then we have $\des(\pi)=2$, and $\ides(\pi)=3$ since $\pi^{-1}=415362$.
It is easy to check that $\ides(\pi)$ also is the number of instances in which ${i}$ appears to the right of ${i+1}$ in $\pi$.

For $1\leq i,j\leq n$, denote by
$$\mathcal{A}_{n, i, j}=\{\pi\in\mathfrak{S}_n\mid \des(\pi)=i-1\text{ and }\ides(\pi)=j-1\},$$
and let $A_{n, i, j}=|\mathcal{A}_{n, i, j}|$.
The generating function of $A_{n,i,j}$,
$$
A_{n}(s, t)=\sum_{\pi \in \mathfrak{S}_{n}} s^{\des\left(\pi\right)+1} t^{\ides(\pi)+1}=\sum_{i, j=1}^{n} A_{n, i, j} s^{i} t^{j},
$$
is called \emph{two-sided Eulerian polynomial} by Petersen \cite{Petersen-2013}
and describes the joint distribution of descents and idescents over the symmetric group.
Clearly, $A_n(t,1)$ gives the classical Eulerian polynomial
\[A_n(t)=\sum_{\pi\in\mathfrak{S}_n}t^{\des(\pi)+1}.\]

The polynomial $A_n(s, t)$ was first studied by Carlitz, Roselle, and Scoville \cite{Carlitz-1966} in 1966, though rather than descents and idescents, they looked at the
equivalent notion of “jumps” to ascents and “readings” to iascents.
They proved the following recurrence for the coefficients $A_{n,i,j}$,
where one could refer \cite[Eq. (7.8)]{Carlitz-1966} for details
but should note that there is an typo in the last row.

\begin{thm}[{\cite[Eq. (10)]{Petersen-2013}}]\label{recA}
For $n\geq 2$ and $1\leq i,j\leq n$, we have
\begin{equation}\label{recAnAn-1}
\begin{aligned}
n A_{n, i, j}
=&(i j+n-1) A_{n-1, i, j}+(1-n+j(n+1-i)) A_{n-1, i-1, j}\\
&+(1-n+i(n+1-j)) A_{n-1, i, j-1}\\
&+(n-1+(n+1-i)(n+1-j)) A_{n-1, i-1, j-1},
\end{aligned}
\end{equation}
where $A_{1,1,1}=1$ and $A_{n,i,j}=0$ if $i\leq 0$ or $j\leq 0$.
\end{thm}


In this paper,
we utilize the geometric representation of permutations
to give an answer to the question posed by Petersen \cite{Petersen-2013}
whether there is a  visual way to understand the recursive relation \eqref{recAnAn-1}.

\begin{definition}[{\cite[Section 1.5]{Stanley-2012}}]
  For a permutation $\pi=\pi_1\pi_2\ldots \pi_n$, the \emph{permutation grid} $P_\pi$ of $\pi$ is an $n\times n$ grid with the $\pi_i$-th (from the left to the right) square in the $i$-th (from the top to the bottom) row filled in.
\end{definition}

For example,
the grid
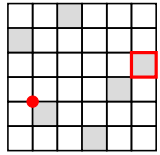
\begin{center}
\begin{tikzpicture}[scale =.65]
\def\hezi{-- +(5mm,0mm) -- +(5mm,5mm) -- +(0mm,5mm) -- cycle [line width=0.6pt]}
\def\judy{-- +(5mm,0mm) -- +(5mm,5mm) -- +(0mm,5mm) -- cycle [line width=0.6pt,fill=gainsboro]}
\tikzstyle{rdot}=[circle,fill=red,draw=red,inner sep=1.5]
\draw (0mm,0mm)\hezi;
\draw (5mm,0mm)\hezi;
\draw (10mm,0mm)\judy;
\draw (15mm,0mm)\hezi;
\draw (20mm,0mm)\hezi;
\draw (25mm,0mm)\hezi;
\draw (0mm,-5mm)\judy;
\draw (5mm,-5mm)\hezi;
\draw (10mm,-5mm)\hezi;
\draw (15mm,-5mm)\hezi;
\draw (20mm,-5mm)\hezi;
\draw (25mm,-5mm)\hezi;
\draw (0mm,-10mm)\hezi;
\draw (5mm,-10mm)\hezi;
\draw (10mm,-10mm)\hezi;
\draw (15mm,-10mm)\hezi;
\draw (20mm,-10mm)\hezi;
\draw (25mm,-10mm)\judy;

\draw (0mm,-15mm)\hezi;
\draw (5mm,-15mm)\hezi;
\draw (10mm,-15mm)\hezi;
\draw (15mm,-15mm)\hezi;
\draw (20mm,-15mm)\judy;
\draw (25mm,-15mm)\hezi;
\draw (0mm,-20mm)\hezi;
\draw (5mm,-20mm)\judy;
\draw (10mm,-20mm)\hezi;
\draw (15mm,-20mm)\hezi;
\draw (20mm,-20mm)\hezi;
\draw (25mm,-20mm)\hezi;
\draw (0mm,-25mm)\hezi;
\draw (5mm,-25mm)\hezi;
\draw (10mm,-25mm)\hezi;
\draw (15mm,-25mm)\judy;
\draw (20mm,-25mm)\hezi;
\draw (25mm,-25mm)\hezi;
\node[rdot] at (5mm, -15mm){};
\draw (25mm,-10mm)-- +(5mm,0mm) -- +(5mm,5mm) -- +(0mm,5mm) -- cycle [line width=1.2pt,red];
\end{tikzpicture}
\end{center}
indicates the permutation $\pi=316524$.
For an $n\times n$ grid,
denote by  $\langle i,j\rangle\in[n]\times[n]$
the square in the $i$-th row and $j$-th column,
and $(i,j)\in[n+1]\times [n+1]$
the grid point intersected by the $i$-th horizontal line and $j$-th vertical line.
Thus, the filled square with the red border in the above grid is indexed by $\langle3,6\rangle$,
and the red grid point is indexed by $(5,2)$.

By investigating geometric forms and propertities
of descents and idescents in permutation grids and
other subtle internal structures,
we give a combinatorial interpretation of the recurrence \eqref{recAnAn-1} in Theorem \ref{recAnAn-1},
which is one of the main results in this paper.

A permutation $\pi\in\mathfrak{S}_n$ is an \emph{involution} if $\pi=\pi^{-1}$,
which implies that the grid $P_\pi$ is symmetrical about the main diagonal.
We call an involution $\sigma\in\mathfrak{S}_n$ is \emph{fixed-point free}
if for any $1\leq i\leq n$ we have $\sigma_i\neq i$.
It is clear that no fixed-point forces
no filled square on the main diagonal of the grid $P_\sigma$.
For example,
the grids of involution $\pi=42513$ and fixed-point free involution $\sigma=532614$ are listed below.
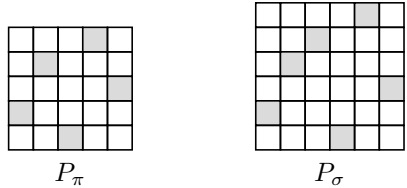
\begin{center}
\begin{tikzpicture}[scale = 0.65]
\def\hezi{-- +(5mm,0mm) -- +(5mm,5mm) -- +(0mm,5mm) -- cycle [line width=0.6pt]}
\def\judy{-- +(5mm,0mm) -- +(5mm,5mm) -- +(0mm,5mm) -- cycle [line width=0.6pt,fill=gainsboro]}
\tikzstyle{rdot}=[circle,fill=red,draw=red,inner sep=1.5]

\draw (0mm,-5mm)\hezi;
\draw (5mm,-5mm)\hezi;
\draw (10mm,-5mm)\hezi;
\draw (15mm,-5mm)\judy;
\draw (20mm,-5mm)\hezi;
\draw (0mm,-10mm)\hezi;
\draw (5mm,-10mm)\judy;
\draw (10mm,-10mm)\hezi;
\draw (15mm,-10mm)\hezi;
\draw (20mm,-10mm)\hezi;
\draw (0mm,-15mm)\hezi;
\draw (5mm,-15mm)\hezi;
\draw (10mm,-15mm)\hezi;
\draw (15mm,-15mm)\hezi;
\draw (20mm,-15mm)\judy;
\draw (0mm,-20mm)\judy;
\draw (5mm,-20mm)\hezi;
\draw (10mm,-20mm)\hezi;
\draw (15mm,-20mm)\hezi;
\draw (20mm,-20mm)\hezi;
\draw (0mm,-25mm)\hezi;
\draw (5mm,-25mm)\hezi;
\draw (10mm,-25mm)\judy;
\draw (15mm,-25mm)\hezi;
\draw (20mm,-25mm)\hezi;
\node at (12.5mm,-30mm){$P_{\pi}$};


\begin{scope}[shift={(50mm,0mm)}]
\draw (0mm,0mm)\hezi;
\draw (5mm,0mm)\hezi;
\draw (10mm,0mm)\hezi;
\draw (15mm,0mm)\hezi;
\draw (20mm,0mm)\judy;
\draw (25mm,0mm)\hezi;
\draw (0mm,-5mm)\hezi;
\draw (5mm,-5mm)\hezi;
\draw (10mm,-5mm)\judy;
\draw (15mm,-5mm)\hezi;
\draw (20mm,-5mm)\hezi;
\draw (25mm,-5mm)\hezi;
\draw (0mm,-10mm)\hezi;
\draw (5mm,-10mm)\judy;
\draw (10mm,-10mm)\hezi;
\draw (15mm,-10mm)\hezi;
\draw (20mm,-10mm)\hezi;
\draw (25mm,-10mm)\hezi;

\draw (0mm,-15mm)\hezi;
\draw (5mm,-15mm)\hezi;
\draw (10mm,-15mm)\hezi;
\draw (15mm,-15mm)\hezi;
\draw (20mm,-15mm)\hezi;
\draw (25mm,-15mm)\judy;
\draw (0mm,-20mm)\judy;
\draw (5mm,-20mm)\hezi;
\draw (10mm,-20mm)\hezi;
\draw (15mm,-20mm)\hezi;
\draw (20mm,-20mm)\hezi;
\draw (25mm,-20mm)\hezi;
\draw (0mm,-25mm)\hezi;
\draw (5mm,-25mm)\hezi;
\draw (10mm,-25mm)\hezi;
\draw (15mm,-25mm)\judy;
\draw (20mm,-25mm)\hezi;
\draw (25mm,-25mm)\hezi;
\node at (15mm,-30mm){$P_{\sigma}$};
\end{scope}
\end{tikzpicture}
\end{center}

For consistency,
we keep the same notations in \cite{Guo-Zeng-2006}.
For $n\geq 1$,
let $\mathcal{I}_{n}$ and $\mathcal{J}_{n}$ be the sets of involutions
and fixed-point free involutions on $[n]$, respectively,
and let
\[
I_{n,k}=\left|\{\pi\in\mathcal{I}_n\mid \des(\pi)=k\}\right|,\qquad
J_{n,k}=\left|\{\pi\in\mathcal{J}_n\mid \des(\pi)=k\}\right|.
\]
Motivated by Brenti's conjecture \cite{Brenti-1994} on the log-concavity of the sequence $I_{n,k}$,
Guo and Zeng \cite{Guo-Zeng-2006} verified that the
two sequences $I_{n,k}$ and $J_{n,k}$ are both unimodal in $k$ for all $n$,
where a sequence $a_0,a_1,\ldots,a_n$ of real numbers is
\emph{log-concave} if $a_i^2\leq a_{i-1}a_{i+1}$ for any $1\leq i\leq n-1$,
and \emph{unimodal} if
 $a_0\leq a_1\leq \cdots \leq a_j\geq a_{j+1}\geq \cdots \geq a_n$ for some $0\leq j\leq n$.
Define
$$
\begin{aligned}
&I_{n}(t)=\sum_{\pi \in \mathcal{I}_{n}} t^{\mathrm{des}(\pi)}=\sum_{k=0}^{n-1} I_{n, k} t^{k}, \\[4pt]
&J_{n}(t)=\sum_{\pi \in \mathcal{J}_{n}} t^{\mathrm{des}(\pi)}=\sum_{k=0}^{n-1} J_{n, k} t^{k}.
\end{aligned}
$$

Based on the generating functions of polynomials $I_n(t)$ and $J_n(t)$,
\begin{align*}
  \sum_{n=0}^{\infty}I_n(t)\frac{u^n}{(1-t)^{n+1}}&=\sum_{r=0}^{\infty}\frac{t^r}{(1-u)^{r+1}(1-u^2)^{r(r+1)/2}},\\[4pt]
  \sum_{n=0}^{\infty}J_n(t)\frac{u^n}{(1-t)^{n+1}}&=\sum_{r=0}^{\infty}\frac{t^r}{(1-u^2)^{r(r+1)/2}},
\end{align*}
obtained by D\'{e}sarm\'{e}nien-Foata \cite{Desarmenien-1985} and Gessel-Reutenauer \cite{Gessel-1993},
Guo and Zeng \cite{Guo-Zeng-2006} derived two linear recurrence formulas for $I_{n,k}$ and $J_{n,k}$ that play a critical role in proving unimodality.

\begin{thm}[\cite{Guo-Zeng-2006}, Theorem 2.2]\label{thmInk}
For $n \geq 3$ and $k \geq 0$, the numbers $I_{n, k}$ satisfy the following recurrence formula
\begin{equation}\label{recInk}
\begin{aligned}
n I_{n, k}=&(k+1) I_{n-1, k}+(n-k) I_{n-1, k-1}+\left((k+1)^{2}+n-2\right) I_{n-2, k} \\[3pt]
&+(2 k(n-k-1)-n+3) I_{n-2, k-1}+\left((n-k)^{2}+n-2\right) I_{n-2, k-2},
\end{aligned}
\end{equation}
where $I_{1,0}=1, I_{2,0}=1, I_{2,1}=1$ and $I_{n, k}=0$ if $k<0$.
\end{thm}

\begin{thm}[\cite{Guo-Zeng-2006}, Theorem 2.1]\label{thmJnk}
For $n \geq 2$ and $k \geq 0$, the numbers $J_{2 n, k}$ satisfy the following recurrence formula
\begin{equation}\label{recJnk}
\begin{aligned}
2 n J_{2 n, k}=& {(k(k+1)+2 n-2) J_{2 n-2, k}+2((k-1)(2 n-k-1)+1) J_{2 n-2, k-1} } \\[3pt]
&+((2 n-k)(2 n-k+1)+2 n-2) J_{2 n-2, k-2},
\end{aligned}
\end{equation}
where $J_{2,0}=0, J_{2,1}=1$ and $J_{2 n, k}=0$ if $k<0$.
\end{thm}

They stated that
it would be interesting to find a combinatorial proof of the recurrence formulas \eqref{recInk} and \eqref{recJnk},
since such a proof could hopefully lead to a combinatorial proof of the unimodality of these two
sequences.
Notice that the idescents of an involution are coincident with the descents.
Thus by developing the constructive analysis that arises from the proof of Theorem \ref{recAnAn-1},
we prove Theorems \ref{thmInk} and \ref{thmJnk} combinatorially.



The rest of this paper is organized as follows.
Section \ref{secpg} is served as preparation for introducing
internal structures of permutation grids and related operations on grids.
In Section \ref{secA},
we give the combinatorial proof of Theorem \ref{recA}.
In Section \ref{Secinvolu}, building on the work of the previous section,
we prove Theorems \ref{thmInk} and \ref{thmJnk} combinatorially.

\section{Internal structures of permutation grids}\label{secpg}

In this section, we first introduce the inserting operations and corresponding deleting operations on permutation grids.
With the inserting operations,
we observe the changes of descents and idescents.



For any $\pi\in\mathfrak{S}_n$ and $1\leq i,j\leq n+1$,
the \emph{inserting operation} $\varphi_{(i,j)}$
is defined by inserting a filled square at the grid point $(i,j)$ of grid $P_\pi$,
while keep the relative positions of original filled squares.


\begin{exam}
Let $\pi=361542$, then $\varphi_{(3,5)}(\pi)=\sigma=3751642$.
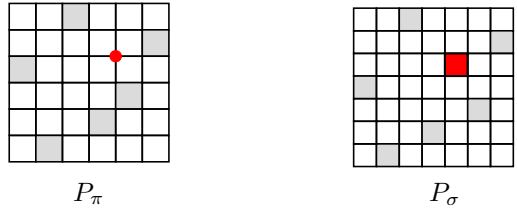
\begin{center}
\begin{tikzpicture}[scale = 0.7]
\def\hezi{-- +(5mm,0mm) -- +(5mm,5mm) -- +(0mm,5mm) -- cycle [line width=0.6pt]}
\def\judy{-- +(5mm,0mm) -- +(5mm,5mm) -- +(0mm,5mm) -- cycle [line width=0.6pt,fill=gainsboro]}
\tikzstyle{rdot}=[circle,fill=red,draw=red,inner sep=1.5]
\tikzstyle{gdot}=[circle,fill=babyblue,draw=babyblue,inner sep=1.5]

\draw (0mm,0mm)\hezi;
\draw (5mm,0mm)\hezi;
\draw (10mm,0mm)\judy;
\draw (15mm,0mm)\hezi;
\draw (20mm,0mm)\hezi;
\draw (25mm,0mm)\hezi;
\draw (0mm,-5mm)\hezi;
\draw (5mm,-5mm)\hezi;
\draw (10mm,-5mm)\hezi;
\draw (15mm,-5mm)\hezi;
\draw (20mm,-5mm)\hezi;
\draw (25mm,-5mm)\judy;
\draw (0mm,-10mm)\judy;
\draw (5mm,-10mm)\hezi;
\draw (10mm,-10mm)\hezi;
\draw (15mm,-10mm)\hezi;
\draw (20mm,-10mm)\hezi;
\draw (25mm,-10mm)\hezi;
\draw (0mm,-15mm)\hezi;
\draw (5mm,-15mm)\hezi;
\draw (10mm,-15mm)\hezi;
\draw (15mm,-15mm)\hezi;
\draw (20mm,-15mm)\judy;
\draw (25mm,-15mm)\hezi;
\draw (0mm,-20mm)\hezi;
\draw (5mm,-20mm)\hezi;
\draw (10mm,-20mm)\hezi;
\draw (15mm,-20mm)\judy;
\draw (20mm,-20mm)\hezi;
\draw (25mm,-20mm)\hezi;
\draw (0mm,-25mm)\hezi;
\draw (5mm,-25mm)\judy;
\draw (10mm,-25mm)\hezi;
\draw (15mm,-25mm)\hezi;
\draw (20mm,-25mm)\hezi;
\draw (25mm,-25mm)\hezi;
\node[rdot] at (20mm, -5mm){};
\node at (15mm,-31mm){$P_{\pi}$};
\end{tikzpicture}
\hspace{6em}
\begin{tikzpicture}[scale = 0.6]
\def\hezi{-- +(5mm,0mm) -- +(5mm,5mm) -- +(0mm,5mm) -- cycle [line width=0.6pt]}
\def\judy{-- +(5mm,0mm) -- +(5mm,5mm) -- +(0mm,5mm) -- cycle [line width=0.6pt,fill=gainsboro]}
\tikzstyle{rdot}=[circle,fill=red,draw=red,inner sep=1.5]
\tikzstyle{gdot}=[circle,fill=babyblue,draw=babyblue,inner sep=1.5]

\draw (0mm,0mm)\hezi;
\draw (5mm,0mm)\hezi;
\draw (10mm,0mm)\judy;
\draw (15mm,0mm)\hezi;
\draw (20mm,0mm)\hezi;
\draw (25mm,0mm)\hezi;
\draw (30mm,0mm)\hezi;
\draw (0mm,-5mm)\hezi;
\draw (5mm,-5mm)\hezi;
\draw (10mm,-5mm)\hezi;
\draw (15mm,-5mm)\hezi;
\draw (20mm,-5mm)\hezi;
\draw (25mm,-5mm)\hezi;
\draw (30mm,-5mm)\judy;
\draw (0mm,-10mm)\hezi;
\draw (5mm,-10mm)\hezi;
\draw (10mm,-10mm)\hezi;
\draw (15mm,-10mm)\hezi;
\draw (20mm,-10mm)-- +(5mm,0mm) -- +(5mm,5mm) -- +(0mm,5mm) -- cycle [line width=0.8pt,fill=red];
\draw (25mm,-10mm)\hezi;
\draw (30mm,-10mm)\hezi;
\draw (0mm,-15mm)\judy;
\draw (5mm,-15mm)\hezi;
\draw (10mm,-15mm)\hezi;
\draw (15mm,-15mm)\hezi;
\draw (20mm,-15mm)\hezi;
\draw (25mm,-15mm)\hezi;
\draw (30mm,-15mm)\hezi;
\draw (0mm,-20mm)\hezi;
\draw (5mm,-20mm)\hezi;
\draw (10mm,-20mm)\hezi;
\draw (15mm,-20mm)\hezi;
\draw (20mm,-20mm)\hezi;
\draw (25mm,-20mm)\judy;
\draw (30mm,-20mm)\hezi;
\draw (0mm,-25mm)\hezi;
\draw (5mm,-25mm)\hezi;
\draw (10mm,-25mm)\hezi;
\draw (15mm,-25mm)\judy;
\draw (20mm,-25mm)\hezi;
\draw (25mm,-25mm)\hezi;
\draw (30mm,-25mm)\hezi;
\draw (0mm,-30mm)\hezi;
\draw (5mm,-30mm)\judy;
\draw (10mm,-30mm)\hezi;
\draw (15mm,-30mm)\hezi;
\draw (20mm,-30mm)\hezi;
\draw (25mm,-30mm)\hezi;
\draw (30mm,-30mm)\hezi;

\node at (20mm,-36mm){$P_\sigma$};
\end{tikzpicture}

\end{center}
\end{exam}

Clearly, the above insertion is roughly reversible.
For any $\sigma\in \mathfrak{S}_{n+1}$ with filled square $\ij$ in the grid $P_\sigma$,
we can delete the $i$-th row and the $j$-th column from $P_\sigma$
to produce a new permutation grid $P_\pi$ with $\pi\in\mathfrak{S}_{n}$.
We call this operation as \emph{deleting operation},
which can be regarded as the inverse operation of $\varphi_{(i,j)}$
and denoted by $\varphi^{-1}_{\ij}$.
Hence, we have $\pi=\varphi^{-1}_{\ij}(\sigma)$
if and only if $\sigma=\varphi_{(i,j)}(\pi)$.


Note that in the grid of $\pi=\pi_1\pi_2\cdots \pi_n$,
each descent $i$ appears
when the filled square $\langle i,\pi_i\rangle $ is located to the right of the filled square $\langle i+1,\pi_{i+1}\rangle$.
Hence we call two filled squares of the above relative position a \emph{des-pair}, otherwise an \emph{asc-pair}.
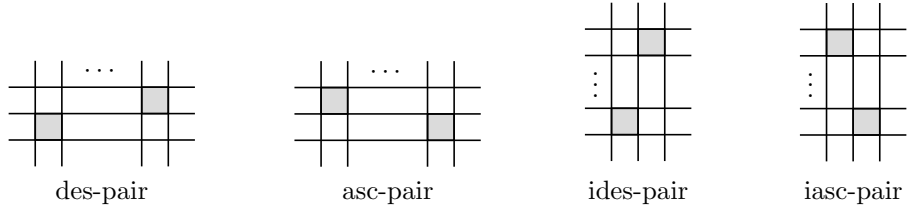
\begin{center}
\begin{tikzpicture}[scale = 0.7]
\def\hezi{-- +(5mm,0mm) -- +(5mm,5mm) -- +(0mm,5mm) -- cycle [line width=0.6pt]}
\def\judy{-- +(5mm,0mm) -- +(5mm,5mm) -- +(0mm,5mm) -- cycle [line width=0.6pt,fill=gainsboro]}
\tikzstyle{bdot}=[circle,fill=black,draw=black,inner sep=1.2]
\tikzstyle{gdot}=[circle,fill=babyblue,draw=babyblue,inner sep=1.5]
\draw (5mm,-5mm)\judy;
\draw (25mm,0mm)\judy;
\draw[line width=0.6pt] (0mm,5mm)--(35mm,5mm);
\draw[line width=0.6pt] (0mm,0mm)--(35mm,0mm);
\draw[line width=0.6pt] (0mm,-5mm)--(35mm,-5mm);
\draw[line width=0.6pt] (5mm,-10mm)--(5mm,10mm);
\draw[line width=0.6pt] (10mm,-10mm)--(10mm,10mm);
\draw[line width=0.6pt] (25mm,-10mm)--(25mm,10mm);
\draw[line width=0.6pt] (30mm,-10mm)--(30mm,10mm);
\node at (17.5mm,8mm){$\cdots$};
\node at (17.5mm,-15mm){des-pair};
\end{tikzpicture}
\hspace{3em}
\begin{tikzpicture}[scale = 0.7]
\def\hezi{-- +(5mm,0mm) -- +(5mm,5mm) -- +(0mm,5mm) -- cycle [line width=0.6pt]}
\def\judy{-- +(5mm,0mm) -- +(5mm,5mm) -- +(0mm,5mm) -- cycle [line width=0.6pt,fill=gainsboro]}
\tikzstyle{bdot}=[circle,fill=black,draw=black,inner sep=1.2]
\tikzstyle{gdot}=[circle,fill=babyblue,draw=babyblue,inner sep=1.5]
\draw (5mm,0mm)\judy;
\draw (25mm,-5mm)\judy;
\draw[line width=0.6pt] (0mm,5mm)--(35mm,5mm);
\draw[line width=0.6pt] (0mm,0mm)--(35mm,0mm);
\draw[line width=0.6pt] (0mm,-5mm)--(35mm,-5mm);
\draw[line width=0.6pt] (5mm,-10mm)--(5mm,10mm);
\draw[line width=0.6pt] (10mm,-10mm)--(10mm,10mm);
\draw[line width=0.6pt] (25mm,-10mm)--(25mm,10mm);
\draw[line width=0.6pt] (30mm,-10mm)--(30mm,10mm);
\node at (17.5mm,8mm){$\cdots$};
\node at (17.5mm,-15mm){asc-pair};
\end{tikzpicture}
\hspace{3em}
\begin{tikzpicture}[scale = 0.7]
\def\hezi{-- +(5mm,0mm) -- +(5mm,5mm) -- +(0mm,5mm) -- cycle [line width=0.6pt]}
\def\judy{-- +(5mm,0mm) -- +(5mm,5mm) -- +(0mm,5mm) -- cycle [line width=0.6pt,fill=gainsboro]}
\tikzstyle{bdot}=[circle,fill=black,draw=black,inner sep=1.2]
\tikzstyle{gdot}=[circle,fill=babyblue,draw=babyblue,inner sep=1.5]
\draw (0mm,0mm)\judy;
\draw (5mm,15mm)\judy;
\draw[line width=0.6pt] (5mm,-5mm)--(5mm,25mm);
\draw[line width=0.6pt] (0mm,-5mm)--(0mm,25mm);
\draw[line width=0.6pt] (10mm,-5mm)--(10mm,25mm);
\draw[line width=0.6pt] (-5mm,20mm)--(15mm,20mm);
\draw[line width=0.6pt] (-5mm,15mm)--(15mm,15mm);
\draw[line width=0.6pt] (-5mm,5mm)--(15mm,5mm);
\draw[line width=0.6pt] (-5mm,0mm)--(15mm,0mm);
\node at (-3mm,11mm){$\vdots$};
\node at (5mm,-11mm){ides-pair};
\end{tikzpicture}
\hspace{3em}
\begin{tikzpicture}[scale = 0.7]
\def\hezi{-- +(5mm,0mm) -- +(5mm,5mm) -- +(0mm,5mm) -- cycle [line width=0.6pt]}
\def\judy{-- +(5mm,0mm) -- +(5mm,5mm) -- +(0mm,5mm) -- cycle [line width=0.6pt,fill=gainsboro]}
\tikzstyle{bdot}=[circle,fill=black,draw=black,inner sep=1.2]
\tikzstyle{gdot}=[circle,fill=babyblue,draw=babyblue,inner sep=1.5]
\draw (5mm,0mm)\judy;
\draw (0mm,15mm)\judy;
\draw[line width=0.6pt] (5mm,-5mm)--(5mm,25mm);
\draw[line width=0.6pt] (0mm,-5mm)--(0mm,25mm);
\draw[line width=0.6pt] (10mm,-5mm)--(10mm,25mm);
\draw[line width=0.6pt] (-5mm,20mm)--(15mm,20mm);
\draw[line width=0.6pt] (-5mm,15mm)--(15mm,15mm);
\draw[line width=0.6pt] (-5mm,5mm)--(15mm,5mm);
\draw[line width=0.6pt] (-5mm,0mm)--(15mm,0mm);
\node at (-3mm,11mm){$\vdots$};
\node at (5mm,-11mm){iasc-pair};
\end{tikzpicture}
\end{center}
Similarly, an idescent $j$ appears in  $P_\pi$ if the filled square $\langle \pi^{-1}_j ,j\rangle$ is
below the filled square $\langle \pi^{-1}_{j+1},j+1\rangle$.
We call these two squares an \emph{ides-pair},
otherwise an \emph{iasc-pair}.

Particularly, we call a des-pair that is also a ides-pair  a $21$-\emph{pair},
and an asc-pair that is also an iasc-pair a $12$-\emph{pair}.
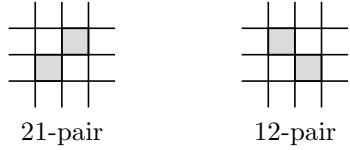
\begin{center}
\begin{tikzpicture}[scale = 0.7]
\def\hezi{-- +(5mm,0mm) -- +(5mm,5mm) -- +(0mm,5mm) -- cycle [line width=0.6pt]}
\def\judy{-- +(5mm,0mm) -- +(5mm,5mm) -- +(0mm,5mm) -- cycle [line width=0.6pt,fill=gainsboro]}
\tikzstyle{bdot}=[circle,fill=black,draw=black,inner sep=1.2]
\tikzstyle{gdot}=[circle,fill=babyblue,draw=babyblue,inner sep=1.5]
\draw (5mm,-5mm)\judy;
\draw (10mm,0mm)\judy;
\draw[line width=0.6pt] (0mm,5mm)--(20mm,5mm);
\draw[line width=0.6pt] (0mm,0mm)--(20mm,0mm);
\draw[line width=0.6pt] (0mm,-5mm)--(20mm,-5mm);
\draw[line width=0.6pt] (5mm,-10mm)--(5mm,10mm);
\draw[line width=0.6pt] (15mm,-10mm)--(15mm,10mm);
\draw[line width=0.6pt] (10mm,-10mm)--(10mm,10mm);
\node at (10mm,-15mm){$21$-pair};
\end{tikzpicture}
\hspace{4em}
\begin{tikzpicture}[scale = 0.7]
\def\hezi{-- +(5mm,0mm) -- +(5mm,5mm) -- +(0mm,5mm) -- cycle [line width=0.6pt]}
\def\judy{-- +(5mm,0mm) -- +(5mm,5mm) -- +(0mm,5mm) -- cycle [line width=0.6pt,fill=gainsboro]}
\tikzstyle{bdot}=[circle,fill=black,draw=black,inner sep=1.2]
\tikzstyle{gdot}=[circle,fill=babyblue,draw=babyblue,inner sep=1.5]
\draw (5mm,0mm)\judy;
\draw (10mm,-5mm)\judy;
\draw[line width=0.6pt] (0mm,5mm)--(20mm,5mm);
\draw[line width=0.6pt] (0mm,0mm)--(20mm,0mm);
\draw[line width=0.6pt] (0mm,-5mm)--(20mm,-5mm);
\draw[line width=0.6pt] (5mm,-10mm)--(5mm,10mm);
\draw[line width=0.6pt] (15mm,-10mm)--(15mm,10mm);
\draw[line width=0.6pt] (10mm,-10mm)--(10mm,10mm);
\node at (10mm,-15mm){$12$-pair};
\end{tikzpicture}
\end{center}


\begin{definition}\label{defdtype}
  Let $\pi\in\mathfrak{S}_n$ and $P_\pi$ be its grid,
  for $1\leq i,j\leq n+1$,
  we say the grid point $(i,j)$ is of \emph{$d$-type} $(p,q)$ if
  $$\left(\des\left(\sigma\right),\ides\left(\sigma\right)\right)-(\des(\pi),\ides(\pi))=(p,q),$$
  where $\sigma=\varphi_{(i,j)}(\pi)$.
  Moreover, we call $p$ the \emph{$d_h$-type} and $q$ the \emph{$d_v$-type}.
\end{definition}

Obviously, the $d$-type records the changes of numbers of descents and idescents
after the insertion at the grid point $(i,j)$.
The following propositions assert that the distributions of grid points
of certain $d_h$-types and $d_v$-types are affected by des-pairs, ides-pairs, asc-pairs and iasc-paris.

\begin{prop}\label{prop-dhdes}
The $d_h$-types of the grid points on the middle horizontal  line of
des-pairs and asc-pairs are indicated in the figure below,
where the numbers $0$ and $1$ indicate the $d_h$-types of the grid points
on the segments enclosed by the corresponding braces.
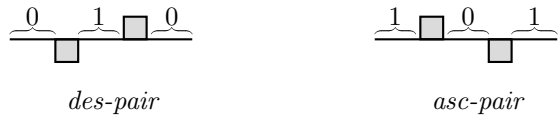
\begin{center}
\begin{tikzpicture}[scale =0.6]
\def\judy{-- +(5mm,0mm) -- +(5mm,5mm) -- +(0mm,5mm) -- cycle [line width=0.8pt, fill=gainsboro]}
\def\nbx{-- +(5mm,0mm) -- +(5mm,5mm) -- +(0mm,5mm) -- cycle [line width=0.6pt,pattern={Lines[angle=45,distance=2.5pt, line width=0.6pt]}, pattern color=black]}
\def\bbrace{[decorate,pen colour={black},
    decoration = {calligraphic brace,mirror,raise=0.7mm,aspect=0.5}]}


\draw(0mm,0mm)\judy;
\draw(15mm,5mm)\judy;
\draw[line width=0.8pt](-10mm,5mm)--(30mm,5mm);


\draw\bbrace  (0mm,5mm) --  (-10mm,5mm);
\node at (-5mm,10mm){$0$};

\draw\bbrace (14mm,5mm) --  (5mm,5mm);
\node at (9.5mm,10mm){$1$};

\draw\bbrace (30mm,5mm) --  (21mm,5mm);
\node at (25.5mm,10mm){$0$};

\node at (12.5mm,-9mm){des-pair};


\begin{scope}[shift={(80mm,0mm)}]
\draw(0mm,5mm)\judy;
\draw(15mm,0mm)\judy;
\draw[line width=0.8pt](-10mm,5mm)--(30mm,5mm);


\draw\bbrace (-1mm,5mm) --  (-10mm,5mm);
\node at (-5.5mm,10mm){$1$};

\draw\bbrace (15mm,5mm) --  (6mm,5mm);
\node at (10.5mm,10mm){$0$};

\draw\bbrace (30mm,5mm) --  (20mm,5mm);
\node at (25mm,10mm){$1$};

\node at (12.5mm,-9mm){asc-pair};
\end{scope}

\end{tikzpicture}
\end{center}

\end{prop}

\begin{prop}\label{prop-dvdes}
The $d_v$-types of the grid points on the middle vertical  line of
ides-pairs and iasc-pairs are indicated in the figure below,
where the numbers $0$ and $1$ indicate the $d_v$-types of the grid points on the segments enclosed by the corresponding braces.
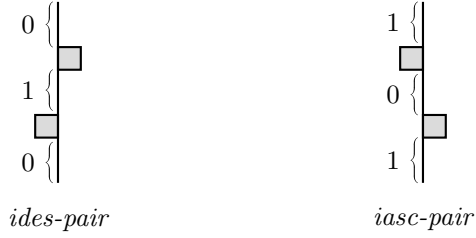
\begin{center}
\begin{tikzpicture}[scale =0.6]
\def\judy{-- +(5mm,0mm) -- +(5mm,5mm) -- +(0mm,5mm) -- cycle [line width=0.8pt, fill=gainsboro]}
\def\nbx{-- +(5mm,0mm) -- +(5mm,5mm) -- +(0mm,5mm) -- cycle [line width=0.6pt,pattern={Lines[angle=45,distance=2.5pt, line width=0.6pt]}, pattern color=black]}
\def\bbrace{[decorate,pen colour={black},
    decoration = {calligraphic brace,mirror,raise=0.7mm,aspect=0.5}]}
\def\rbrace{[decorate,pen colour={red},
    decoration = {calligraphic brace,mirror,raise=0.7mm,aspect=0.5}]}


\draw(10mm,10mm)\judy;
\draw(5mm,-5mm)\judy;
\draw[line width=0.8pt](10mm,-15mm)--(10mm,25mm);


\draw\bbrace (10mm,25mm) --  (10mm,15mm);
\node at (3.5mm,20mm){$0$};

\draw\bbrace (10mm,10mm) --  (10mm,1mm);
\node at (3.5mm,5.5mm){$1$};

\draw\bbrace (10mm,-6mm) --  (10mm,-15mm);
\node at (3.5mm,-10.5mm){$0$};

\node at (10mm,-23mm){ides-pair};


\begin{scope}[shift={(80mm,0mm)}]
\draw(5mm,10mm)\judy;
\draw(10mm,-5mm)\judy;
\draw[line width=0.8pt](10mm,-15mm)--(10mm,25mm);


\draw\bbrace (10mm,25mm) --  (10mm,16mm);
\node at (3.5mm,20.5mm){$1$};

\draw\bbrace (10mm,9mm) --  (10mm,0mm);
\node at (3.5mm,4.5mm){$0$};

\draw\bbrace (10mm,-5mm) --  (10mm,-15mm);
\node at (3.5mm,-10mm){$1$};

\node at (10mm,-23mm){iasc-pair};
\end{scope}

\end{tikzpicture}
\end{center}

\end{prop}

These two propositions can be directly checked by applying inserting operation $\varphi_{(i,j)}$.
By the same analysis, we characterize the distributions of the $d$-types of the grid points
on the grid boundary.

\begin{prop}\label{prop-dtype-rim}
The $d_h$-types of the grid points on the top and bottom boundaries of the grid,
and the $d_v$-types of the grid points on the left and right boundaries of the grid
are indicated in the figure below, respectively.
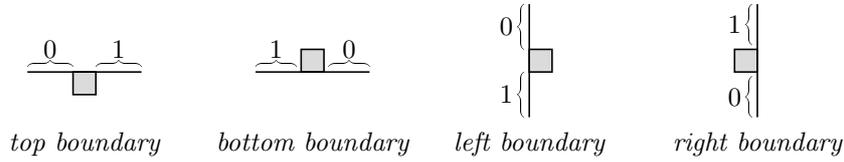
\begin{center}
\begin{tikzpicture}[scale =0.6]
\def\nbx{-- +(5mm,0mm) -- +(5mm,5mm) -- +(0mm,5mm) -- cycle [line width=0.6pt,pattern={Lines[angle=45,distance=2.5pt, line width=0.6pt]}, pattern color=black]}
\def\judy{-- +(5mm,0mm) -- +(5mm,5mm) -- +(0mm,5mm) -- cycle [line width=0.6pt,fill=gainsboro]}
\def\bbrace{[decorate,pen colour={black}, decoration = {calligraphic brace,mirror,raise=0.7mm,aspect=0.5}]}
\def\rbrace{[decorate,pen colour={red}, decoration = {calligraphic brace,mirror,raise=0.7mm,aspect=0.5}]}

\draw (10mm,0mm)\judy;
\draw[line width=0.6pt] (0mm,5mm)--(25mm,5mm);

\draw\bbrace (10mm,5mm) --  (0mm,5mm);
\node at (5mm,10mm){$0$};
\draw\bbrace (25mm,5mm) --  (15mm,5mm);
\node at (20mm,10mm){$1$};

\node at (12.5mm,-11mm){top boundary};
\begin{scope}[xshift=50mm]
\draw (10mm,5mm)\judy;
\draw[line width=0.6pt] (0mm,5mm)--(25mm,5mm);

\draw\bbrace (9mm,5mm) --  (0mm,5mm);
\node at (4.5mm,10mm){$1$};
\draw\bbrace (25mm,5mm) --  (16mm,5mm);
\node at (20.5mm,10mm){$0$};

\node at (12.5mm,-11mm){bottom boundary};
\end{scope}

\begin{scope}[xshift=100mm]
\draw (10mm,5mm)\judy;
\draw[line width=0.6pt] (10mm,-5mm)--(10mm,20mm);

\draw\bbrace (10mm,20mm) --  (10mm,10mm);
\node at (5mm,15mm){$0$};
\draw\bbrace (10mm,5mm) --  (10mm,-5mm);
\node at (5mm,0mm){$1$};


\node at (10mm,-11mm){left boundary};
\end{scope}

\begin{scope}[xshift=150mm]

\draw (5mm,5mm)\judy;
\draw[line width=0.6pt] (10mm,-5mm)--(10mm,20mm);

\draw\bbrace (10mm,20mm) --  (10mm,11mm);
\node at (5mm,15.5mm){$1$};
\draw\bbrace (10mm,4mm) --  (10mm,-5mm);
\node at (5mm,-0.5mm){$0$};

\node at (10mm,-11mm){right boundary};

\end{scope}
\end{tikzpicture}
\end{center}
\end{prop}

Based on the above observations,
for any $\pi\in\mathfrak{S}_n$,
and any grid point $(i,j)$ in $P_\pi$ for $1\leq i,j\leq n+1$,
we have
\begin{equation}\label{0and1}
(\mathrm{des}(\sigma),\mathrm{ides}(\sigma))-(\mathrm{des}(\pi),\mathrm{ides}(\pi))\in\{(0,0), (0,1), (1,0), (1,1)\},
\end{equation}
where $\sigma=\varphi_{(i,j)}(\pi)$.
Or equivalently, in any permutation grid,
all possible $d$-types of grid points only are
$(0,0)$, $(1,0),(0,1)$ and $(1,1)$.


In a permutation grid, for $p\in\{0,1\}$ (resp. $q\in\{0,1\}$),
we first join the grid points of $d_h$-type $p$ (resp. $d_v$-type $q$) only by horizontal (resp. vertical) lines,
and then connect such two grid points once they are on the same filled square,
which creates several paths containing all possible grid points with $d_h$-type $p$ (resp. $d_v$-type $q$),
and we call theses paths by \emph{$p_h$-paths} (resp. \emph{$q_v$-paths}).
It follows from \eqref{0and1} that each grid point must belong to
only one $p_h$-path, and only one $q_v$-path.

\begin{exam}
For $\pi=316524$, we present the grid $P_\pi$ with
each one of $0_h$-paths, $1_h$-paths, $0_v$-paths and $1_v$-paths as follows.
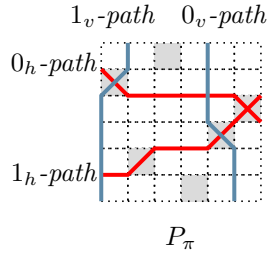
\begin{center}
\begin{tikzpicture}[scale =0.7]
\def\hezi{-- +(5mm,0mm) -- +(5mm,5mm) -- +(0mm,5mm) -- cycle [line width=0.6pt, dotted]}
\def\judy{-- +(5mm,0mm) -- +(5mm,5mm) -- +(0mm,5mm) -- cycle [line width=0.6pt,fill=gainsboro,dotted]}
\tikzstyle{cc}=[circle,draw=black,fill=yellow, line width=0.5pt, inner sep=1.5]
\draw (0mm,0mm)\hezi;
\draw (5mm,0mm)\hezi;
\draw (10mm,0mm)\judy;
\draw (15mm,0mm)\hezi;
\draw (20mm,0mm)\hezi;
\draw (25mm,0mm)\hezi;
\draw (0mm,-5mm)\judy;
\draw (5mm,-5mm)\hezi;
\draw (10mm,-5mm)\hezi;
\draw (15mm,-5mm)\hezi;
\draw (20mm,-5mm)\hezi;
\draw (25mm,-5mm)\hezi;
\draw (0mm,-10mm)\hezi;
\draw (5mm,-10mm)\hezi;
\draw (10mm,-10mm)\hezi;
\draw (15mm,-10mm)\hezi;
\draw (20mm,-10mm)\hezi;
\draw (25mm,-10mm)\judy;
\draw (0mm,-15mm)\hezi;
\draw (5mm,-15mm)\hezi;
\draw (10mm,-15mm)\hezi;
\draw (15mm,-15mm)\hezi;
\draw (20mm,-15mm)\judy;
\draw (25mm,-15mm)\hezi;
\draw (0mm,-20mm)\hezi;
\draw (5mm,-20mm)\judy;
\draw (10mm,-20mm)\hezi;
\draw (15mm,-20mm)\hezi;
\draw (20mm,-20mm)\hezi;
\draw (25mm,-20mm)\hezi;
\draw (0mm,-25mm)\hezi;
\draw (5mm,-25mm)\hezi;
\draw (10mm,-25mm)\hezi;
\draw (15mm,-25mm)\judy;
\draw (20mm,-25mm)\hezi;
\draw (25mm,-25mm)\hezi;
\node at (-9mm,1mm){$0_h$-path};
\node at (-9mm,-20mm){$1_h$-path};
\draw(0mm,0mm)--(5mm,-5mm)--(25mm,-5mm)--(30mm,-10mm)[line width=1.5pt,draw=red];
\draw (0mm,-20mm)--(5mm,-20mm)--(10mm,-15mm)--(20mm,-15mm)--(25mm,-10mm)--(30mm,-5mm)[line width=1.5pt,draw=red];
\node at (2mm,10mm){$1_v$-path};
\node at (23mm,10mm){$0_v$-path};
\draw(20mm,5mm)--(20mm,-10mm)--(25mm,-15mm)--(25mm,-25mm)[line width=1.5pt,draw=afblue];
\draw (5mm,5mm)--(5mm,0mm)--(0mm,-5mm)--(0mm,-25mm)[line width=1.5pt,draw=afblue];

\node at (15mm,-32mm){$P_\pi$};
\end{tikzpicture}
\end{center}
\end{exam}

The next lemma determines the specific walking trends of the $p_h$-paths and $q_v$-paths.

\begin{lem}\label{lempath}
In any permutation grid $P_\pi$,
\begin{itemize}
\item
each $0_h$ (resp. $1_h$)-path goes from the left boundary to the right  boundary along the horizontal grid lines
except for carrying out a southeast (resp. northeast) step when encountering a filled square;
\item
each $0_v$ (resp. $1_v$)-path goes from the top boundary to the bottom boundary along the vertical grid lines
except for carrying out a southeast (resp. southwest) step when encountering a filled square.
\end{itemize}
\end{lem}

\pf
Let $\pi\in\mathfrak{S}_n$ and $P_\pi$ be its grid, for $1\leq i,j\leq n$,
assume that the square $\ij$ is filled.
We only need to determine the distributions of $d_h$-types and $d_v$-types
in the neighborhood of the square $\ij$,
where the neighborhood contains grid points
on the segments that directly touch the four corners of the square $\ij$,
and it is clear that these segments are formed
either by the square $\ij$ and the filled squares on its adjacent rows or columns,
or by the square $\ij$ and the grid boundary.
See the following figure for a better understanding.

\begin{center}
\begin{tikzpicture}[scale =0.8]
\def\hezi{-- +(5mm,0mm) -- +(5mm,5mm) -- +(0mm,5mm) -- cycle [line width=0.8pt, densely dotted, fill=gainsboro]}
\def\judy{-- +(5mm,0mm) -- +(5mm,5mm) -- +(0mm,5mm) -- cycle [line width=0.8pt,fill=gainsboro]}
\tikzstyle{rdot}=[circle,fill=red,draw=red,inner sep=0pt, outer sep=0pt,minimum size=1mm]
\tikzstyle{gdot}=[circle,fill=afblue,draw=afblue,inner sep=0pt, outer sep=0pt, ,minimum size=1mm]
\tikzstyle{halfg}=[semicircle,fill=afblue,draw=afblue,inner sep=0pt, outer sep=0pt, ,minimum size=0.5mm,anchor=south]

\draw(15mm,0mm)\judy;
\draw[line width=0.8pt](-3mm,5mm)--(38mm,5mm);
\draw[line width=0.8pt](-3mm,0mm)--(38mm,0mm);

\draw(0mm,5mm)\hezi;\draw(0mm,-5mm)\hezi;\draw(30mm,5mm)\hezi;\draw(30mm,-5mm)\hezi;

\draw[line width=0.8pt] (15mm,-18mm)--(15mm,23mm);
\draw[line width=0.8pt] (20mm,-18mm)--(20mm,23mm);

\draw(10mm,-15mm)\hezi;\draw(20mm,-15mm)\hezi;\draw(10mm,15mm)\hezi;\draw(20mm,15mm)\hezi;

\node(s1)[rdot] at (15mm,0mm){};\node[rdot] at (10mm,0mm){};\node[rdot] at (5mm,0mm){};
\draw [decorate,pen colour={red},
    decoration = {calligraphic brace,mirror,raise=0.7mm,aspect=0.25}] (5mm,0mm) --  (14mm,0mm);
\node at (7.5mm,-4mm){\red{$1$}};

\node[rdot] at (20mm,5mm){};\node[rdot] at (25mm,5mm){};\node[rdot] at (30mm,5mm){};
\draw [decorate,pen colour={red},
    decoration = {calligraphic brace,mirror,raise=0.7mm,aspect=0.25}] (30mm,5mm) --  (21mm,5mm);
\node at (27.5mm,9mm){\red{$1$}};


\node[rdot] at (15mm,5mm){};\node[rdot] at (10mm,5mm){};\node[rdot] at (5mm,5mm){};
\draw [decorate,pen colour={red},
    decoration = {calligraphic brace,mirror,raise=0.7mm,aspect=0.75}] (14mm,5mm) --  (5mm,5mm);
\node at (7.5mm,9mm){\red{$0$}};

\node[rdot] at (20mm,0mm){};\node[rdot] at (25mm,0mm){};\node[rdot] at (30mm,0mm){};
\draw [decorate,pen colour={red},
    decoration = {calligraphic brace,mirror,raise=0.7mm,aspect=0.75}] (21mm,0mm) --  (30mm,0mm);
\node at (27.5mm,-4mm){\red{$0$}};

\node[halfg,at={(s1.center)},rotate=-135] at (15mm,0mm){};\node[gdot] at (15mm,-5mm){};\node[gdot] at (15mm,-10mm){};

\draw [decorate,pen colour={afblue},
    decoration = {calligraphic brace,mirror,raise=0.7mm,aspect=0.75}] (15mm,-1mm) --  (15mm,-10mm);
\node at (11.5mm,-7.5mm){\afblue{$1$}};

\node[halfg,at={(s1.center)},rotate=45] at (20mm,5mm){};\node[gdot] at (20mm,10mm){};\node[gdot] at (20mm,15mm){};

\draw [decorate,pen colour={afblue},
    decoration = {calligraphic brace,mirror,raise=0.7mm,aspect=0.75}] (20mm,6mm) --  (20mm,15mm);
\node at (23.5mm,12.5mm){\afblue{$1$}};

\node[halfg,at={(s1.center)},rotate=135] at (20mm,0mm){};\node[gdot] at (20mm,-5mm){};\node[gdot] at (20mm,-10mm){};

\draw [decorate,pen colour={afblue},
    decoration = {calligraphic brace,mirror,raise=0.7mm,aspect=0.25}] (20mm,-10mm) --  (20mm,-1mm);
\node at (23.5mm,-7.5mm){\afblue{$0$}};

\node[halfg,at={(s1.center)},rotate=-45] at (15mm,5mm){};\node[gdot] at (15mm,10mm){};\node[gdot] at (15mm,15mm){};

\draw [decorate,pen colour={afblue},
    decoration = {calligraphic brace,mirror,raise=0.7mm,aspect=0.25}] (15mm,15mm) --  (15mm,6mm);
\node at (11.5mm,12.5mm){\afblue{$0$}};


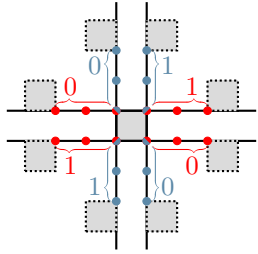
\node at (17.5mm,-23mm){the filled square $\ij$ with its neighborhood };

\end{tikzpicture}
\end{center}

Note that the square $\ij$ together with the filled squares in its adjacent rows or columns must
composes one certain of des-pair, asc-pair, ides-pair, and iasc-pair.
Thus combining Propositions \ref{prop-dhdes}--\ref{prop-dtype-rim},
we derive the distributions of $d_h$-types and $d_v$-types
as labeled in the figure,
where the $d_h$-types are labeled in red and the $d_v$-types in blue.
%
%
%
Therefore, according to the constructions of $p_h$-paths and $q_v$-paths for $p,q\in\{0,1\}$,
we evidently completes the proof.
\qed

\begin{remark}\label{rem0011}
In the proof of Lemma \ref{lempath},
we also see that for any permutation grid $P_\pi$ and any filled square $\ij$,
the grid points $(i,j)$ and $(i+1,j+1)$ always receive the $d$-type $(0,0)$,
and the grid points $(i+1,j)$ and $(i,j+1)$ always receive the $d$-type $(1,1)$.
The above fact can be also observed by using the structures of $12$-pairs and $21$-pairs when analysing the behaviors of inserting filled squares at the four corners of the filled square $\ij$.

\end{remark}

According to the specific walking trends of the $p_h$-paths and $q_v$-paths for $p,q\in\{0,1\}$ in the grid $P_\pi$,
we can use $\des(\pi)$ and $\ides(\pi)$ to calculate the numbers of these paths.

\begin{thm}\label{thmpaths}
Let $\pi\in\mathfrak{S}_n$, then in the  grid $P_\pi$,
  \begin{enumerate}[(a)]
    \item the number of $0_h$-paths is counted by $\des(\pi)+1$,
    \item the number of $1_h$-paths is counted by $n-\des(\pi)$,
    \item the number of $0_v$-paths is counted by $\ides(\pi)+1$,
    \item the number of $1_v$-paths is counted by $n-\ides(\pi)$.
  \end{enumerate}
\end{thm}

\proof
With the help of Lemma \ref{lempath},
we know that for $p\in\{0,1\}$ (resp. $q\in\{0,1\}$),
the number of $p_h$-paths (resp. $q_v$-paths) is determined by
the number of grid points of $d_h$-type $p$ (resp. $d_v$-type $q$) in the first vertical (resp. horizontal) line of the grid $P_\pi$.

By Proposition \ref{prop-dhdes},
for $2\leq i\leq n$,
the number of grid points $(i,1)$ of $d_h$-type $0$
is equal to $\des(\pi)$, the number of des-pairs,
and the number of $(i,1)$ of $d_h$-type $1$ is equal to $(n-1)-\des(\pi)$,
the number of asc-pairs.
And by Proposition \ref{prop-dtype-rim}, the grid points $(1,1)$ and $(n+1,1)$
always receive $d_h$-type $0$ and $1$, respectively.
Therefore, the number of $0_h$-paths is counted by $\des(\pi)+1$,
and the number of $1_h$-path is counted by $(n-1)-\des(\pi)+1$.
Using the similar analysis, by Propositions \ref{prop-dvdes} and \ref{prop-dtype-rim},
we can verify the statements $(c)$ and $(d)$.\qed


\section{Recurrence of $A_{n,i,j}$}\label{secA}

In this section, we first compute the numbers of grid points of certain $d$-types in the permutation grid with given numbers of descents and idescents,
then give a combinatorial proof of the recurrence for $A_{n,i,j}$
in Theorem \ref{recA}.
Recall that for $1\leq i,j\leq n$,
$$
\mathcal{A}_{n, i, j}=\left\{\pi \in \mathfrak{S}_{n}: \des\left(\pi\right)=i-1 \text { and } \ides(\pi)=j-1\right\} .
$$

\begin{thm}\label{thmdtype}
  For any $\pi\in \mathcal{A}_{n,i,j}$, in the permutation grid $P_\pi$, there are
  \begin{enumerate}[(a)]
    \item $(ij+n)$ grid points of d-type $(0,0)$,\label{crossa}
    \item $(j(n+1-i)-n)$ grid points of d-type $(1,0)$,
    \item $(i(n+1-j)-n)$ grid points of d-type $(0,1)$,
    \item $((n+1-i)(n+1-j)+n)$ grid points of d-type $(1,1)$.
  \end{enumerate}
\end{thm}

\pf
It is evident that for $p,q\in \{0,1\}$,
a grid point is of $d$-type $(p,q)$
if and only if it is the intersecting point of a $p_h$-path and a $q_v$-path.
Since $\pi\in \mathcal{A}_{n,i,j}$,
we have $\des(\pi)=i-1$ and $\ides(\pi)=j-1$,
and there are $i$ $0_h$-paths, $j$ $0_v$-paths, $(n+1-i)$ $1_h$-paths and $(n+1-j)$ $1_v$-paths in $P_\pi$ by Theorem \ref{thmpaths}.

In the following proof,
one can get a clearer picture from Example \ref{examdtype}.
Since $i$ $0_h$-paths travel from the left to the right
and $j$ $0_v$-paths travel from the top to the bottom,
these paths should encounter each other $ij$ times.
However, by Lemma \ref{lempath},
the $0_h$-paths and the $0_v$-paths both take one northeast step once touching a filled square,
which implies that they meet only once in each filled square, but form two intersecting crosses on the grid points.
Therefore, the $i$ $0_h$-paths and the $j$ $0_v$-paths yield $(ij+n)$ grid points of $d$-type $(0,0)$, which verifies the statement $(a)$.
And the statement $(d)$ can be proved similarly by replacing $0_h$-paths and $0_v$-paths with $1_h$-paths and $1_v$-paths, respectively.

To see the number of grid points of $d$-type $(1,0)$, with Lemma \ref{lempath},
note that when encountering filled squares,
the $1_h$-path carries out a northeast step,
while the $0_v$-path carries out a southeast step,
which leads to that the intersecting points of
$1_h$-paths and $0_v$-paths caused by filled squares
are inside the squares but not on the grid points.
Hence, $(n+1-i)$ $1_h$-paths and $j$ $0_v$-paths cross $j(n+1-i)$ times,
but only give $(j(n+1-i)-n)$ grid points of $d$-type $(1,0)$,
which proves the statements $(b)$.
Similarly, $i$ $0_h$-paths and $(n+1-j)$ $1_v$-paths cross $i(n+1-j)$ times but only form $(i(n+1-j)-n)$ grid points of $d$-type $(0,1)$,
which gives the statement $(c)$.\qed


\begin{exam}\label{examdtype}
Let $\pi=316524\in\mathcal{A}_{6,4,4}$, the following grids
show the $d$-type of each grid point in $P_\pi$ by Theorem \ref{thmdtype}.
\begin{center}
\begin{tikzpicture}[scale =0.7]
\def\hezi{-- +(5mm,0mm) -- +(5mm,5mm) -- +(0mm,5mm) -- cycle [line width=0.6pt, dotted]}
\def\judy{-- +(5mm,0mm) -- +(5mm,5mm) -- +(0mm,5mm) -- cycle [line width=0.6pt,fill=gainsboro,dotted]}
\tikzstyle{cc}=[circle,draw=black,fill=yellow, line width=0.5pt, inner sep=1.5]
\draw (0mm,0mm)\hezi;
\draw (5mm,0mm)\hezi;
\draw (10mm,0mm)\judy;
\draw (15mm,0mm)\hezi;
\draw (20mm,0mm)\hezi;
\draw (25mm,0mm)\hezi;
\draw (0mm,-5mm)\judy;
\draw (5mm,-5mm)\hezi;
\draw (10mm,-5mm)\hezi;
\draw (15mm,-5mm)\hezi;
\draw (20mm,-5mm)\hezi;
\draw (25mm,-5mm)\hezi;
\draw (0mm,-10mm)\hezi;
\draw (5mm,-10mm)\hezi;
\draw (10mm,-10mm)\hezi;
\draw (15mm,-10mm)\hezi;
\draw (20mm,-10mm)\hezi;
\draw (25mm,-10mm)\judy;
\draw (0mm,-15mm)\hezi;
\draw (5mm,-15mm)\hezi;
\draw (10mm,-15mm)\hezi;
\draw (15mm,-15mm)\hezi;
\draw (20mm,-15mm)\judy;
\draw (25mm,-15mm)\hezi;
\draw (0mm,-20mm)\hezi;
\draw (5mm,-20mm)\judy;
\draw (10mm,-20mm)\hezi;
\draw (15mm,-20mm)\hezi;
\draw (20mm,-20mm)\hezi;
\draw (25mm,-20mm)\hezi;
\draw (0mm,-25mm)\hezi;
\draw (5mm,-25mm)\hezi;
\draw (10mm,-25mm)\hezi;
\draw (15mm,-25mm)\judy;
\draw (20mm,-25mm)\hezi;
\draw (25mm,-25mm)\hezi;
\node at (-3mm,5mm){\red{$0$}};
\node at (-3mm,0mm){\red{$0$}};
\node at (-3mm,-5mm){\red{$1$}};
\node at (-3mm,-10mm){\red{$0$}};
\node at (-3mm,-15mm){\red{$0$}};
\node at (-3mm,-20mm){\red{$1$}};
\node at (-3mm,-25mm){\red{$1$}};
\draw[shift={(0,-0.05)}] (0mm,5mm)--(10mm,5mm)--(15mm,0mm)--(30mm,0mm)[line width=1.5pt,draw=red];
\draw[shift={(0,-0.05)}] (0mm,0mm)--(5mm,-5mm)--(25mm,-5mm)--(30mm,-10mm)[line width=1.5pt,draw=red];
\draw[shift={(0,-0.05)}] (0mm,-10mm)--(20mm,-10mm)--(25mm,-15mm)--(30mm,-15mm)[line width=1.5pt,draw=red];
\draw[shift={(0,-0.05)}] (0mm,-15mm)--(5mm,-15mm)--(10mm,-20mm)--(15mm,-20mm)--(20mm,-25mm)--(30mm,-25mm)[line width=1.5pt,draw=red];
\node at (0mm,9mm){\afblue{$0$}};
\node at (5mm,9mm){\afblue{$1$}};
\node at (10mm,9mm){\afblue{$0$}};
\node at (15mm,9mm){\afblue{$1$}};
\node at (20mm,9mm){\afblue{$0$}};
\node at (25mm,9mm){\afblue{$0$}};
\node at (30mm,9mm){\afblue{$1$}};
\draw[shift={(0.05,0)}] (0mm,5mm)--(0mm,0mm)--(5mm,-5mm)--(5mm,-15mm)--(10mm,-20mm)--(10mm,-25mm)[line width=1.5pt,draw=afblue];
\draw[shift={(0.05,0)}] (10mm,5mm)--(15mm,0mm)--(15mm,-20mm)--(20mm,-25mm)[line width=1.5pt,draw=afblue];
\draw[shift={(0.05,0)}] (20mm,5mm)--(20mm,-10mm)--(25mm,-15mm)--(25mm,-25mm)[line width=1.5pt,draw=afblue];
\draw[shift={(0.05,0)}] (25mm,5mm)--(25mm,-5mm)--(30mm,-10mm)--(30mm,-25mm)[line width=1.5pt,draw=afblue];
\node[cc] at (0mm, 5mm){};\node[cc] at (10mm, 5mm){};\node[cc] at (15mm, 0mm){};\node[cc] at (20mm, 0mm){};\node[cc] at (25mm, 0mm){};

\node[cc] at (0mm, 0mm){};\node[cc] at (5mm, -5mm){};\node[cc] at (15mm, -5mm){};\node[cc] at (20mm, -5mm){};\node[cc] at (25mm, -5mm){};\node[cc] at (30mm, -10mm){};

\node[cc] at (5mm, -10mm){};\node[cc] at (15mm, -10mm){};\node[cc] at (20mm, -10mm){};\node[cc] at (25mm, -15mm){};\node[cc] at (30mm, -15mm){};

\node[cc] at (5mm, -15mm){};\node[cc] at (10mm, -20mm){};\node[cc] at (15mm, -20mm){};\node[cc] at (20mm, -25mm){};\node[cc] at (25mm, -25mm){};\node[cc] at (30mm, -25mm){};

\node at (15mm,-32mm){$22$ grid points of $d$-type $(0,0)$};
\end{tikzpicture}
\hspace{3em}
\begin{tikzpicture}[scale =0.7]
\def\hezi{-- +(5mm,0mm) -- +(5mm,5mm) -- +(0mm,5mm) -- cycle [line width=0.6pt, dotted]}
\def\judy{-- +(5mm,0mm) -- +(5mm,5mm) -- +(0mm,5mm) -- cycle [line width=0.6pt,fill=gainsboro,dotted]}
\tikzstyle{cc}=[circle,draw=black,fill=yellow, line width=0.5pt, inner sep=1.5]
\draw (0mm,0mm)\hezi;
\draw (5mm,0mm)\hezi;
\draw (10mm,0mm)\judy;
\draw (15mm,0mm)\hezi;
\draw (20mm,0mm)\hezi;
\draw (25mm,0mm)\hezi;
\draw (0mm,-5mm)\judy;
\draw (5mm,-5mm)\hezi;
\draw (10mm,-5mm)\hezi;
\draw (15mm,-5mm)\hezi;
\draw (20mm,-5mm)\hezi;
\draw (25mm,-5mm)\hezi;
\draw (0mm,-10mm)\hezi;
\draw (5mm,-10mm)\hezi;
\draw (10mm,-10mm)\hezi;
\draw (15mm,-10mm)\hezi;
\draw (20mm,-10mm)\hezi;
\draw (25mm,-10mm)\judy;
\draw (0mm,-15mm)\hezi;
\draw (5mm,-15mm)\hezi;
\draw (10mm,-15mm)\hezi;
\draw (15mm,-15mm)\hezi;
\draw (20mm,-15mm)\judy;
\draw (25mm,-15mm)\hezi;
\draw (0mm,-20mm)\hezi;
\draw (5mm,-20mm)\judy;
\draw (10mm,-20mm)\hezi;
\draw (15mm,-20mm)\hezi;
\draw (20mm,-20mm)\hezi;
\draw (25mm,-20mm)\hezi;
\draw (0mm,-25mm)\hezi;
\draw (5mm,-25mm)\hezi;
\draw (10mm,-25mm)\hezi;
\draw (15mm,-25mm)\judy;
\draw (20mm,-25mm)\hezi;
\draw (25mm,-25mm)\hezi;
\node at (-3mm,5mm){\red{$0$}};
\node at (-3mm,0mm){\red{$0$}};
\node at (-3mm,-5mm){\red{$1$}};
\node at (-3mm,-10mm){\red{$0$}};
\node at (-3mm,-15mm){\red{$0$}};
\node at (-3mm,-20mm){\red{$1$}};
\node at (-3mm,-25mm){\red{$1$}};
\draw (0mm,-5mm)--(5mm,0mm)--(10mm,0mm)--(15mm,5mm)--(30mm,5mm)[line width=1.5pt,draw=red];
\draw (0mm,-20mm)--(5mm,-20mm)--(10mm,-15mm)--(20mm,-15mm)--(25mm,-10mm)--(30mm,-5mm)[line width=1.5pt,draw=red];
\draw(0mm,-25mm)--(15mm,-25mm)--(20mm,-20mm)--(30mm,-20mm)[line width=1.5pt,draw=red];
\node at (0mm,9mm){\afblue{$0$}};
\node at (5mm,9mm){\afblue{$1$}};
\node at (10mm,9mm){\afblue{$0$}};
\node at (15mm,9mm){\afblue{$1$}};
\node at (20mm,9mm){\afblue{$0$}};
\node at (25mm,9mm){\afblue{$0$}};
\node at (30mm,9mm){\afblue{$1$}};
\draw (0mm,5mm)--(0mm,0mm)--(5mm,-5mm)--(5mm,-15mm)--(10mm,-20mm)--(10mm,-25mm)[line width=1.5pt,draw=afblue];
\draw (10mm,5mm)--(15mm,0mm)--(15mm,-20mm)--(20mm,-25mm)[line width=1.5pt,draw=afblue];
\draw(20mm,5mm)--(20mm,-10mm)--(25mm,-15mm)--(25mm,-25mm)[line width=1.5pt,draw=afblue];
\draw(25mm,5mm)--(25mm,-5mm)--(30mm,-10mm)--(30mm,-25mm)[line width=1.5pt,draw=afblue];
\node[cc] at (20mm, 5mm){};\node[cc] at (25mm, 5mm){};

\node[cc] at (15mm, -15mm){};

\node[cc] at (10mm, -25mm){};\node[cc] at (25mm, -20mm){};\node[cc] at (30mm, -20mm){};

\node at (15mm,-32mm){$6$ grid points of $d$-type $(1,0)$};
\end{tikzpicture}
\end{center}

\begin{center}
\begin{tikzpicture}[scale =0.7]
\def\hezi{-- +(5mm,0mm) -- +(5mm,5mm) -- +(0mm,5mm) -- cycle [line width=0.6pt, dotted]}
\def\judy{-- +(5mm,0mm) -- +(5mm,5mm) -- +(0mm,5mm) -- cycle [line width=0.6pt,fill=gainsboro,dotted]}
\tikzstyle{cc}=[circle,draw=black,fill=yellow, line width=0.5pt, inner sep=1.5]
\draw (0mm,0mm)\hezi;
\draw (5mm,0mm)\hezi;
\draw (10mm,0mm)\judy;
\draw (15mm,0mm)\hezi;
\draw (20mm,0mm)\hezi;
\draw (25mm,0mm)\hezi;
\draw (0mm,-5mm)\judy;
\draw (5mm,-5mm)\hezi;
\draw (10mm,-5mm)\hezi;
\draw (15mm,-5mm)\hezi;
\draw (20mm,-5mm)\hezi;
\draw (25mm,-5mm)\hezi;
\draw (0mm,-10mm)\hezi;
\draw (5mm,-10mm)\hezi;
\draw (10mm,-10mm)\hezi;
\draw (15mm,-10mm)\hezi;
\draw (20mm,-10mm)\hezi;
\draw (25mm,-10mm)\judy;
\draw (0mm,-15mm)\hezi;
\draw (5mm,-15mm)\hezi;
\draw (10mm,-15mm)\hezi;
\draw (15mm,-15mm)\hezi;
\draw (20mm,-15mm)\judy;
\draw (25mm,-15mm)\hezi;
\draw (0mm,-20mm)\hezi;
\draw (5mm,-20mm)\judy;
\draw (10mm,-20mm)\hezi;
\draw (15mm,-20mm)\hezi;
\draw (20mm,-20mm)\hezi;
\draw (25mm,-20mm)\hezi;
\draw (0mm,-25mm)\hezi;
\draw (5mm,-25mm)\hezi;
\draw (10mm,-25mm)\hezi;
\draw (15mm,-25mm)\judy;
\draw (20mm,-25mm)\hezi;
\draw (25mm,-25mm)\hezi;
\node at (-3mm,5mm){\red{$0$}};
\node at (-3mm,0mm){\red{$0$}};
\node at (-3mm,-5mm){\red{$1$}};
\node at (-3mm,-10mm){\red{$0$}};
\node at (-3mm,-15mm){\red{$0$}};
\node at (-3mm,-20mm){\red{$1$}};
\node at (-3mm,-25mm){\red{$1$}};
\draw (0mm,5mm)--(10mm,5mm)--(15mm,0mm)--(30mm,0mm)[line width=1.5pt,draw=red];
\draw(0mm,0mm)--(5mm,-5mm)--(25mm,-5mm)--(30mm,-10mm)[line width=1.5pt,draw=red];
\draw(0mm,-10mm)--(20mm,-10mm)--(25mm,-15mm)--(30mm,-15mm)[line width=1.5pt,draw=red];
\draw(0mm,-15mm)--(5mm,-15mm)--(10mm,-20mm)--(15mm,-20mm)--(20mm,-25mm)--(30mm,-25mm)[line width=1.5pt,draw=red];
\node at (0mm,9mm){\afblue{$0$}};
\node at (5mm,9mm){\afblue{$1$}};
\node at (10mm,9mm){\afblue{$0$}};
\node at (15mm,9mm){\afblue{$1$}};
\node at (20mm,9mm){\afblue{$0$}};
\node at (25mm,9mm){\afblue{$0$}};
\node at (30mm,9mm){\afblue{$1$}};
\draw (5mm,5mm)--(5mm,0mm)--(0mm,-5mm)--(0mm,-25mm)[line width=1.5pt,draw=afblue];
\draw (15mm,5mm)--(10mm,0mm)--(10mm,-15mm)--(5mm,-20mm)--(5mm,-25mm)[line width=1.5pt,draw=afblue];
\draw (30mm,5mm)--(30mm,-5mm)--(20mm,-15mm)--(20mm,-20mm)--(15mm,-25mm)[line width=1.5pt,draw=afblue];

\node[cc] at (5mm, 5mm){};\node[cc] at (30mm, 0mm){};

\node[cc] at (10mm, -5mm){};

\node[cc] at (0mm, -10mm){};\node[cc] at (10mm, -10mm){};

\node[cc] at (0mm, -15mm){};\node[cc] at (10mm, -10mm){};

\node at (15mm,-32mm){$6$ grid points of $d$-type $(0,1)$};
\end{tikzpicture}
\hspace{3em}
\begin{tikzpicture}[scale =0.7]
\def\hezi{-- +(5mm,0mm) -- +(5mm,5mm) -- +(0mm,5mm) -- cycle [line width=0.6pt, dotted]}
\def\judy{-- +(5mm,0mm) -- +(5mm,5mm) -- +(0mm,5mm) -- cycle [line width=0.6pt,fill=gainsboro,dotted]}
\tikzstyle{cc}=[circle,draw=black,fill=yellow, line width=0.5pt, inner sep=1.5]
\draw (0mm,0mm)\hezi;
\draw (5mm,0mm)\hezi;
\draw (10mm,0mm)\judy;
\draw (15mm,0mm)\hezi;
\draw (20mm,0mm)\hezi;
\draw (25mm,0mm)\hezi;
\draw (0mm,-5mm)\judy;
\draw (5mm,-5mm)\hezi;
\draw (10mm,-5mm)\hezi;
\draw (15mm,-5mm)\hezi;
\draw (20mm,-5mm)\hezi;
\draw (25mm,-5mm)\hezi;
\draw (0mm,-10mm)\hezi;
\draw (5mm,-10mm)\hezi;
\draw (10mm,-10mm)\hezi;
\draw (15mm,-10mm)\hezi;
\draw (20mm,-10mm)\hezi;
\draw (25mm,-10mm)\judy;
\draw (0mm,-15mm)\hezi;
\draw (5mm,-15mm)\hezi;
\draw (10mm,-15mm)\hezi;
\draw (15mm,-15mm)\hezi;
\draw (20mm,-15mm)\judy;
\draw (25mm,-15mm)\hezi;
\draw (0mm,-20mm)\hezi;
\draw (5mm,-20mm)\judy;
\draw (10mm,-20mm)\hezi;
\draw (15mm,-20mm)\hezi;
\draw (20mm,-20mm)\hezi;
\draw (25mm,-20mm)\hezi;
\draw (0mm,-25mm)\hezi;
\draw (5mm,-25mm)\hezi;
\draw (10mm,-25mm)\hezi;
\draw (15mm,-25mm)\judy;
\draw (20mm,-25mm)\hezi;
\draw (25mm,-25mm)\hezi;
\node at (-3mm,5mm){\red{$0$}};
\node at (-3mm,0mm){\red{$0$}};
\node at (-3mm,-5mm){\red{$1$}};
\node at (-3mm,-10mm){\red{$0$}};
\node at (-3mm,-15mm){\red{$0$}};
\node at (-3mm,-20mm){\red{$1$}};
\node at (-3mm,-25mm){\red{$1$}};
\draw[shift={(0,-0.05)}] (0mm,-5mm)--(5mm,0mm)--(10mm,0mm)--(15mm,5mm)--(30mm,5mm)[line width=1.5pt,draw=red];
\draw[shift={(0,-0.05)}] (0mm,-20mm)--(5mm,-20mm)--(10mm,-15mm)--(20mm,-15mm)--(25mm,-10mm)--(30mm,-5mm)[line width=1.5pt,draw=red];
\draw[shift={(0,-0.05)}](0mm,-25mm)--(15mm,-25mm)--(20mm,-20mm)--(30mm,-20mm)[line width=1.5pt,draw=red];
\node at (0mm,9mm){\afblue{$0$}};
\node at (5mm,9mm){\afblue{$1$}};
\node at (10mm,9mm){\afblue{$0$}};
\node at (15mm,9mm){\afblue{$1$}};
\node at (20mm,9mm){\afblue{$0$}};
\node at (25mm,9mm){\afblue{$0$}};
\node at (30mm,9mm){\afblue{$1$}};
\draw[shift={(-0.05,0)}] (5mm,5mm)--(5mm,0mm)--(0mm,-5mm)--(0mm,-25mm)[line width=1.5pt,draw=afblue];
\draw[shift={(-0.05,0)}] (15mm,5mm)--(10mm,0mm)--(10mm,-15mm)--(5mm,-20mm)--(5mm,-25mm)[line width=1.5pt,draw=afblue];
\draw[shift={(-0.05,0)}] (30mm,5mm)--(30mm,-5mm)--(20mm,-15mm)--(20mm,-20mm)--(15mm,-25mm)[line width=1.5pt,draw=afblue];

\node[cc] at (0mm, -5mm){};\node[cc] at (5mm, 0mm){};\node[cc] at (10mm, 0mm){};\node[cc] at (15mm, 5mm){};\node[cc] at (30mm, 5mm){};

\node[cc] at (0mm, -20mm){};\node[cc] at (5mm, -20mm){};\node[cc] at (10mm, -15mm){};\node[cc] at (20mm, -15mm){};\node[cc] at (25mm, -10mm){};\node[cc] at (30mm, -5mm){};

\node[cc] at (0mm, -25mm){};\node[cc] at (5mm, -25mm){};\node[cc] at (15mm, -25mm){};\node[cc] at (20mm, -20mm){};

\node at (15mm,-32mm){$15$ grid points of $d$-type $(1,1)$};
\end{tikzpicture}
\end{center}
\end{exam}

By the definition of $d$-type and Theorem \ref{thmdtype}, for any $\pi\in\mathfrak{S}_n$,
we finally determine the changes of the numbers of $\des(\pi)$ and $\ides(\pi)$
after employing operations $\varphi_{(r,s)}$ on the grid $P_\pi$.

\begin{thm}\label{thmAntoAn+1}
  Let $1\leq i,j\leq n$ and $\pi\in\mathcal{A}_{n,i,j}$, then after employing the operations $\varphi_{(r,s)}$ on $P_\pi$ for $1\leq r,s\leq n+1$,
  we have
 \begin{enumerate}[(a)]
   \item $(ij+n)$ permutations in $\mathcal{A}_{n+1,i,j}$, 
   \item $(j(n+1-i)-n)$ permutations in $\mathcal{A}_{n+1,i+1,j}$, 
   \item  $(i(n+1-j)-n)$ permutations in $\mathcal{A}_{n+1,i,j+1}$, 
   \item $((n+1-i)(n+1-j)+n)$ permutations in $\mathcal{A}_{n+1,i+1,j+1}$. 
 \end{enumerate}
\end{thm}


We now give a combinatorial interpretation for the recurrence \eqref{recA} in Theorem \ref{recA}.

{\noindent{\emph{Combinatorial Proof of Theorem \ref{recA}.}\hskip 2pt}}
Let
$$
\widetilde{\mathcal{A}}_{n,i,j}=\{(\pi, k) \mid \pi\in\mathcal{A}_{n, i, j} \text{ and }1\leq k\leq n \}.
$$
It is trivial that $|\widetilde{\mathcal{A}}_{n,i,j}|=nA_{n,i,j}$,
which is the left side of \eqref{recA}.

To establish the combinatorial relation between the two sides of $\eqref{recAnAn-1}$,
for $1\leq i,j\leq n$ and $p,q\in\{0,1\}$, we define
\begin{equation*}
\mathcal{B}_{n,i,j}^{(p,q)}=\left\{(\sigma,(r,s))\mid \text{ $\sigma\in\mathcal{A}_{n,i,j}$ and $(r,s)$ is  the grid point of $d$-type $(p,q)$ in $P_\sigma$}\right\}.
\end{equation*}
Then by Theorem \ref{thmdtype}, we have
\begin{equation*}
  \begin{aligned}
    &\left|\mathcal{B}_{n-1,i,j}^{(0,0)}\right|=(ij+n-1)A_{n-1,i,j},\\[3pt]
    &\left|\mathcal{B}_{n-1,i-1,j}^{(1,0)}\right|=(j(n-i+1)-n+1)A_{n-1,i-1,j},\\[3pt]
    &\left|\mathcal{B}_{n-1,i,j-1}^{(0,1)}\right|=(i(n-j+1)-n+1)A_{n-1,i,j-1},\\[3pt]
    &\left|\mathcal{B}_{n-1,i-1,j-1}^{(1,1)}\right|=((n+1-i)(n+1-j)-n+1)A_{n-1,i-1,j-1},
  \end{aligned}
\end{equation*}
which correspond to the four terms on the right side of \eqref{recA}.
The above four sets are apparently pairwise disjoint.
Using $\uplus$ represents the union of two disjoint sets,
we write
\[
\biguplus_{p,q\in\{0,1\}}\mathcal{B}_{n-1,i-p,j-q}^{(p,q)}
=\mathcal{B}_{n-1,i,j}^{(0,0)}\uplus \mathcal{B}_{n-1,i-1,j}^{(1,0)}\uplus \mathcal{B}_{n-1,i,j-1}^{(0,1)}\uplus\mathcal{B}_{n-1,i-1,j-1}^{(1,1)}.
\]

For any pair $(\pi,k)\in\widetilde{\mathcal{A}}_{n,i,j}$,
define the mapping $\Psi$ by
\begin{equation*}
\Psi((\pi,k))=(\sigma, (k,\pi_k)),
\end{equation*}
where  $\sigma=\varphi^{-1}_{\langle k, \pi_k\rangle}(\pi)$.
Since $\sigma=\varphi^{-1}_{\langle k, \pi_k\rangle}(\pi)$ if and only if
$\pi=\varphi_{(k,\pi_k)}(\sigma)$,
by Theorem \ref{thmAntoAn+1},
the permutation $\sigma$ generates $\pi$ with the inserting operations $\varphi_{(k,\pi_k)}$ if and only if
$\sigma\in {\mathcal{A}}_{n-1,i-p,j-q}$
and $(k,\pi_k)$ is of $d$-type $(p,q)$ in $P_\sigma$
for  $p,q\in \{0,1\}$.
Thus
\[
\Psi:\,\widetilde{\mathcal{A}}_{n,i,j}\,\leftrightarrow \biguplus_{p,q\in\{0,1\}}\mathcal{B}_{n-1,i-p,j-q}^{(p,q)}
\]
is indeed a bijection, which completes the proof.\qed

\section{Recurrences of $I_{n,k}$ and $J_{n,k}$ on involutions}\label{Secinvolu}

Recall that $\mathcal{I}_n$ is the set of all involutions on $[n]$.
Since the symmetry of  grids of involutions can
be inherited by the paths and the $d$-types, we have the following proposition.

\begin{prop}\label{lemij=ji}
  Let $\pi\in\mathcal{I}_n$, then for $1\leq i,j \leq n+1$ and $p,q\in\{0,1\}$, in the grid $P_\pi$,
  \begin{enumerate}[(a)]
    \item the grid point $(i,j)$ is of $d$-type $(p,q)$ if and only if
    the grid point $(j,i)$ is of $d$-type   $(q,p)$,
    \item the $d$-types of the grid points on the main diagonal of $P_\pi$ are either $(0,0)$ or $(1,1)$.
  \end{enumerate}
\end{prop}

The following grid give an example to the statement $(a)$,
and the statement $(b)$ is directly from the statement $(a)$ by setting $i=j$.
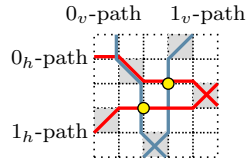
\begin{center}
\begin{tikzpicture}[scale = 0.65]
\def\hezi{-- +(5mm,0mm) -- +(5mm,5mm) -- +(0mm,5mm) -- cycle [line width=0.6pt,dotted]}
\def\judy{-- +(5mm,0mm) -- +(5mm,5mm) -- +(0mm,5mm) -- cycle [line width=0.6pt,fill=gainsboro,dotted]}
\tikzstyle{cc}=[circle,draw=black,fill=yellow, line width=0.5pt, inner sep=1.5]
\tikzstyle{rdot}=[circle,fill=red,draw=red,inner sep=1.5]

\draw (0mm,0mm)\hezi;
\draw (5mm,0mm)\hezi;
\draw (10mm,0mm)\hezi;
\draw (15mm,0mm)\judy;
\draw (20mm,0mm)\hezi;
\draw (0mm,-5mm)\hezi;
\draw (5mm,-5mm)\judy;
\draw (10mm,-5mm)\hezi;
\draw (15mm,-5mm)\hezi;
\draw (20mm,-5mm)\hezi;
\draw (0mm,-10mm)\hezi;
\draw (5mm,-10mm)\hezi;
\draw (10mm,-10mm)\hezi;
\draw (15mm,-10mm)\hezi;
\draw (20mm,-10mm)\judy;
\draw (0mm,-15mm)\judy;
\draw (5mm,-15mm)\hezi;
\draw (10mm,-15mm)\hezi;
\draw (15mm,-15mm)\hezi;
\draw (20mm,-15mm)\hezi;
\draw (0mm,-20mm)\hezi;
\draw (5mm,-20mm)\hezi;
\draw (10mm,-20mm)\judy;
\draw (15mm,-20mm)\hezi;
\draw (20mm,-20mm)\hezi;

\draw[shift={(0,0.05)}] (0mm,0mm)--(5mm,0mm)--(10mm,-5mm)--(20mm,-5mm)--(25mm,-10mm)[line width=1.25pt,draw=red];

\draw (20mm,5mm)--(15mm,0mm)--(15mm,-15mm)--(10mm,-20mm)[line width=1.25pt,draw=afblue];

\node[cc] at (15mm, -5mm){};

\draw[shift={(-0.05,0)}] (5mm,5mm)--(5mm,0mm)--(10mm,-5mm)--(10mm,-15mm)--(15mm,-20mm)[line width=1.25pt,draw=afblue];

\draw (0mm,-15mm)--(5mm,-10mm)--(20mm,-10mm)--(25mm,-5mm)[line width=1.25pt,draw=red];

\node[cc] at (10mm, -10mm){};

\node at (2mm,9mm){\footnotesize $0_v$-path};
\node at (-9mm,0mm){\footnotesize $0_h$-path};

\node at (23mm,9mm){\footnotesize $1_v$-path};
\node at (-9mm,-15mm){\footnotesize $1_h$-path};

\end{tikzpicture}
\end{center}

Hence we can continue to employ inserting operations $\varphi_{(i,j)}$ on some certain grid points in grids of involutions on $[n]$
to generate involutions on $[n+1]$.

\begin{prop}\label{lemiip}
For any $\pi\in\mathcal{I}_n$ and $1\leq i\leq n+1$,
if the grid point $(i,i)$ is of $d$-type $(p,p)$,
then we have
$\varphi_{(i,i)}(\pi)\in \mathcal{I}_{n+1}$,
and
$$\des\left(\varphi_{(i,i)}(\pi)\right)=\des(\pi)+p.$$
\end{prop}

\begin{prop}\label{lemi+1i}
For any $\pi\in\mathcal{I}_n$ and $1\leq i\leq n$,
if the square $\langle i,i\rangle$ is filled in the grid $P_\pi$,
then we have
$\varphi_{(i+1,i)}(\pi)\in \mathcal{I}_{n+1}$, and
$$\des\left(\varphi_{(i+1,i)}(\pi)\right)=\des(\pi)+1.$$
\end{prop}
\proof
Let $\sigma=\varphi_{(i+1,i)}(\pi)$.
The inserting operation $\varphi_{(i+1,i)}$ on $P_\pi$ converts the square $\langle i,i \rangle$
to a symmetric pair of filled squares $\langle i,i+1\rangle$ and $\langle i+1,i\rangle$ in $P_\sigma$.
Thus the grid $P_\sigma$ is still symmetric and $\sigma\in \mathcal{I}_{n+1}$.
Since
the grid point $(i+1,i)$ is of $d$-type $(1,1)$ by Remark \ref{rem0011},
we have $\des\left(\varphi_{(i+1,i)}(\pi)\right)=\des(\pi)+1$.
\qed


We introduce some new operations that insert two filled squares
simultaneously in grids of involutions,
and maintain the symmetry of grids.
Let $\pi\in \mathcal{I}_n$ and $1\leq i,j\leq n+1$ with $i\neq j$,
define \emph{double inserting operation} $\xi_{(i,j)}$
by
inserting filled squares at the grid points $(i,j)$ and $(j,i)$ of the grid $P_\pi$.
Clearly, we have $\xi_{(i,j)}(\pi)\in \mathcal{I}_{n+2}$
and $\xi_{(i,j)}(\pi)=\xi_{(j,i)}(\pi)$.
The following grids show the involution $132$ and the
involution $\xi_{(1,3)}(132)=42513$.
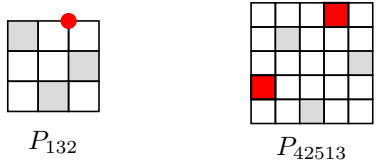
\begin{center}
\begin{tikzpicture}[scale = 0.8]
\def\hezi{-- +(5mm,0mm) -- +(5mm,5mm) -- +(0mm,5mm) -- cycle [line width=0.6pt]}
\def\judy{-- +(5mm,0mm) -- +(5mm,5mm) -- +(0mm,5mm) -- cycle [line width=0.6pt,fill=gainsboro]}

\tikzstyle{rdot}=[circle,fill=red,draw=red,inner sep=2]

\draw (0mm,-15mm)\judy;
\draw (5mm,-15mm)\hezi;
\draw (10mm,-15mm)\hezi;
\draw (0mm,-20mm)\hezi;
\draw (5mm,-20mm)\hezi;
\draw (10mm,-20mm)\judy;
\draw (0mm,-25mm)\hezi;
\draw (5mm,-25mm)\judy;
\draw (10mm,-25mm)\hezi;
\node at (7.5mm,-30mm){$P_{132}$};
\node[rdot] at (10mm,-10mm){};

\begin{scope}[shift={(40mm,-7mm)},scale=.8]
\draw (0mm,-5mm)\hezi;
\draw (5mm,-5mm)\hezi;
\draw (10mm,-5mm)\hezi;
\draw (15mm,-5mm)-- +(5mm,0mm) -- +(5mm,5mm) -- +(0mm,5mm) -- cycle [line width=0.8pt,fill=red];
\draw (20mm,-5mm)\hezi;
\draw (0mm,-10mm)\hezi;
\draw (5mm,-10mm)\judy;
\draw (10mm,-10mm)\hezi;
\draw (15mm,-10mm)\hezi;
\draw (20mm,-10mm)\hezi;
\draw (0mm,-15mm)\hezi;
\draw (5mm,-15mm)\hezi;
\draw (10mm,-15mm)\hezi;
\draw (15mm,-15mm)\hezi;
\draw (20mm,-15mm)\judy;
\draw (0mm,-20mm)-- +(5mm,0mm) -- +(5mm,5mm) -- +(0mm,5mm) -- cycle [line width=0.8pt,fill=red];
\draw (5mm,-20mm)\hezi;
\draw (10mm,-20mm)\hezi;
\draw (15mm,-20mm)\hezi;
\draw (20mm,-20mm)\hezi;
\draw (0mm,-25mm)\hezi;
\draw (5mm,-25mm)\hezi;
\draw (10mm,-25mm)\judy;
\draw (15mm,-25mm)\hezi;
\draw (20mm,-25mm)\hezi;
\node at (12.5mm,-30mm){$P_{42513}$};

%
%
%

\end{scope}
\end{tikzpicture}
\end{center}

For $1<i<n+1$, the \emph{double inserting operations} $\eta_i$ and $\eta'_i$ are defined as
inserting a $12$-pair and a $21$-pair
at the grid point $(i,i)$ in the
permutation grid $P_\pi$, respectively.
Thus $\eta_i(\pi),\eta'_i(\pi)\in\mathcal{I}_{n+2}$ for any $\pi\in\mathcal{I}_n$.
The following example presents the grids of $42315$ and $21354$
that are obtained by employing $\eta_2$ and $\eta'_4$ on
the permutation $213$.
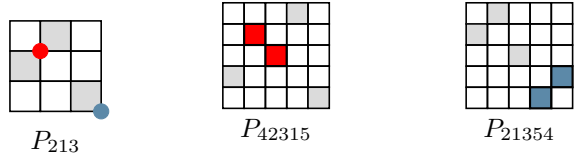
\begin{center}
\begin{tikzpicture}[scale = 0.8]
\def\hezi{-- +(5mm,0mm) -- +(5mm,5mm) -- +(0mm,5mm) -- cycle [line width=0.6pt]}
\def\judy{-- +(5mm,0mm) -- +(5mm,5mm) -- +(0mm,5mm) -- cycle [line width=0.6pt,fill=gainsboro]}

\tikzstyle{rdot}=[circle,fill=red,draw=red,inner sep=2]
\tikzstyle{gdot}=[circle,fill=afblue,draw=afblue,inner sep=2]

\draw (0mm,-15mm)\hezi;
\draw (5mm,-15mm)\judy;
\draw (10mm,-15mm)\hezi;
\draw (0mm,-20mm)\judy;
\draw (5mm,-20mm)\hezi;
\draw (10mm,-20mm)\hezi;
\draw (0mm,-25mm)\hezi;
\draw (5mm,-25mm)\hezi;
\draw (10mm,-25mm)\judy;
\node at (7.5mm,-30mm){$P_{213}$};
\node[rdot] at (5mm,-15mm){};
\node[gdot] at (15mm,-25mm){};

\begin{scope}[shift={(35mm,-7mm)},scale=.7]
\draw (0mm,-5mm)\hezi;
\draw (5mm,-5mm)\hezi;
\draw (10mm,-5mm)\hezi;
\draw (15mm,-5mm)\judy;
\draw (20mm,-5mm)\hezi;
\draw (0mm,-10mm)\hezi;
\draw (5mm,-10mm)-- +(5mm,0mm) -- +(5mm,5mm) -- +(0mm,5mm) -- cycle [line width=0.8pt,fill=red];
\draw (10mm,-10mm)\hezi;
\draw (15mm,-10mm)\hezi;
\draw (20mm,-10mm)\hezi;
\draw (0mm,-15mm)\hezi;
\draw (5mm,-15mm)\hezi;
\draw (10mm,-15mm)-- +(5mm,0mm) -- +(5mm,5mm) -- +(0mm,5mm) -- cycle [line width=0.8pt,fill=red];
\draw (15mm,-15mm)\hezi;
\draw (20mm,-15mm)\hezi;
\draw (0mm,-20mm)\judy;
\draw (5mm,-20mm)\hezi;
\draw (10mm,-20mm)\hezi;
\draw (15mm,-20mm)\hezi;
\draw (20mm,-20mm)\hezi;
\draw (0mm,-25mm)\hezi;
\draw (5mm,-25mm)\hezi;
\draw (10mm,-25mm)\hezi;
\draw (15mm,-25mm)\hezi;
\draw (20mm,-25mm)\judy;
\node at (12.5mm,-30mm){$P_{42315}$};

%
%
%

\end{scope}


\begin{scope}[shift={(75mm,-7mm)},scale=.7]
\draw (0mm,-5mm)\hezi;
\draw (5mm,-5mm)\judy;
\draw (10mm,-5mm)\hezi;
\draw (15mm,-5mm)\hezi;
\draw (20mm,-5mm)\hezi;
\draw (0mm,-10mm)\judy;
\draw (5mm,-10mm)\hezi;
\draw (10mm,-10mm)\hezi;
\draw (15mm,-10mm)\hezi;
\draw (20mm,-10mm)\hezi;
\draw (0mm,-15mm)\hezi;
\draw (5mm,-15mm)\hezi;
\draw (10mm,-15mm)\judy;
\draw (15mm,-15mm)\hezi;
\draw (20mm,-15mm)\hezi;
\draw (0mm,-20mm)\hezi;
\draw (5mm,-20mm)\hezi;
\draw (10mm,-20mm)\hezi;
\draw (15mm,-20mm)\hezi;
\draw (20mm,-20mm)-- +(5mm,0mm) -- +(5mm,5mm) -- +(0mm,5mm) -- cycle [line width=0.8pt,fill=afblue];
\draw (0mm,-25mm)\hezi;
\draw (5mm,-25mm)\hezi;
\draw (10mm,-25mm)\hezi;
\draw (15mm,-25mm)-- +(5mm,0mm) -- +(5mm,5mm) -- +(0mm,5mm) -- cycle [line width=0.8pt,fill=afblue];
\draw (20mm,-25mm)\hezi;

\node at (12.5mm,-30mm){$P_{21354}$};

\end{scope}

\end{tikzpicture}
\end{center}

As the deleting operation $\varphi^{-1}_{\ij}$ to the inserting operation $\varphi_{(i,j)}$,
the double inserting operations $\xi_{(i,j)}, \eta_i$ and $\eta'_i$
also have corresponding deleting operations.

For  $\sigma\in \mathcal{I}_{n+2}$ and $1\leq i,j\leq n+1$ with $i\neq j$,  the \emph{double deleting operation} $\xi^{-1}_{\ij}$ is defined as
deleting the filled squares $\langle i,j+1\rangle$ and $\langle j+1,i\rangle$ from $P_\sigma$ if $i<j$,
and deleting the filled squares $\langle i+1,j\rangle$ and $\langle j,i+1\rangle$ from $P_\sigma$ if $i>j$.
Thus $\xi^{-1}_{\ij}(\sigma)\in \mathcal{I}_n$ and $\xi^{-1}_{\ij}(\sigma)=\xi^{-1}_{\langle j,i\rangle}(\sigma)$.
Clearly, we have $\xi^{-1}_{\ij}(\sigma)=\pi$ if $\xi_{(i,j)}(\pi)=\sigma$.
Note that $\xi^{-1}_{\ij}$ cannot delete any $21$-pair
whose central point is on the main diagonal of $P_\sigma$.

For  $1\leq i\leq n+1$,
\emph{the double deleting operations} $\eta^{-1}_i$ and $\eta'^{-1}_i$ are defined 
by deleting the $12$-pair and $21$-pair whose central grid point is $(i+1, i+1)$ in $P_\sigma$, respectively.
We have $\eta^{-1}_i(\sigma),\eta^{-1}_i(\sigma)\in\mathcal{I}_n$, as well
$\pi=\eta^{-1}_i(\sigma)$ if $\sigma=\eta_i(\pi)$
and $\pi=\eta'^{-1}_i(\sigma)$ if $\sigma=\eta'_i(\pi)$.

Since the descents are coincident with the idescents in involutions,
we can characterize the change of the number of descents of involutions
after applying the above double inserting operations.

\begin{prop}\label{lemxi}
For $\pi\in\mathcal{I}_n$ and $1\leq i\neq j\leq n+1$,
if the grid point $(i,j)$ in $P_\pi$ is of $d$-type $(p,q)$ for $p,q\in\{0,1\}$,
then we have
$$\des\left(\xi_{(i,j)}(\pi)\right)=\des(\pi)+p+q.$$
\end{prop}

\proof
By Proposition \ref{lemij=ji}, the symmetric properties of the permutation grid and the operation $\xi_{(i,j)}$ imply that
the change of descents at the $j$-th row after the insertion at the grid point $(j,i)$ can be completely reflected by
the change of idescents at the $j$-th column after the insertion at the grid point $(i,j)$.
Thus,
we have $\des\left(\xi_{(i,j)}(\pi)\right)=\des(\pi)+p+q$.
\qed

\begin{prop}\label{lemetaeta'}
For $\pi\in\mathcal{I}_n$ and $1\leq i\leq n+1$,
if the grid point $(i,i)$ in $P_\pi$ is of $d$-types $(p,p)$ for $p\in \{0,1\}$,
then we have 
$$\des\left(\eta_i(\pi)\right)=\des(\pi)+p \quad\mbox{and}\quad
\des\left(\eta'_i(\pi)\right)=\des(\pi)+p+1.$$
\end{prop}

\proof
The insertion processes of $\eta_i$ (resp. $\eta'_i$) can be taken as
first inserting a filled square at the grid point $(i,i)$ then extending it to a $12$-pair (resp. $21$-pair).
Thus, the insertion of the filled square should produce $p$ new descent pair at the $i$-th row of $P_\pi$, since the grid point $(i,i)$ is of $d$-type $(p,p)$.
However,
after extending the filled square to a $12$ pair or a $21$ pair,
the possible increase of descents can only come from the interior of the pair.
Notice that the $12$-pair contains no descents
while the $21$-pair contains one descent.
Therefore, we derive that $\des\left(\eta_i(\pi)\right)=\des(\pi)+p$,
and $\des\left(\eta'_i(\pi)\right)=\des(\pi)+p+1$.
\qed


Lightened by the ideas used in the proof of the recurrence of $A_{n,i,j}$,
we shall  combinatorially prove the recursive relation \eqref{recInk} about $I_{n,k}$
in Theorem \ref{thmInk}.

{\noindent{\emph{\textbf{Combinatorial Proof of Theorem \ref{thmInk}.}}\hskip 2pt}}
Let $\mathcal{I}_{n,k}$ be the set of all involutions on $[n]$ with $k$ descents and we see $|\mathcal{I}_{n,k}|=I_{n,k}$,
then the cardinality of the set
\[
\widetilde{\mathcal{I}}_{n,k}=\{(\pi,i)\mid \pi\in \mathcal{I}_{n,k}\text{ and }1\leq i\leq n \}
\]
is $nI_{n,k}$, which equals the left side of \eqref{recInk}.
We next need to construct five pairwise disjoint sets that are counted by the five terms on the right side of \eqref{recInk}.

Let
\begin{equation}\label{defB1nk}
  \begin{aligned}
\mathcal{B}^{(1)}_{n,k}=\{(\sigma,(i,i))\mid&\text{ $\sigma\in\mathcal{I}_{n,k}$,  and $(i,i)$ is the last grid point }\\
    &\text{ at which a $0_h$-path touches the main diagonal of $P_\sigma$}\}.
  \end{aligned}
\end{equation}
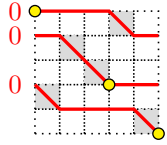
\begin{center}
\begin{tikzpicture}[scale = 0.65]
\def\hezi{-- +(5mm,0mm) -- +(5mm,5mm) -- +(0mm,5mm) -- cycle [line width=0.6pt,dotted]}
\def\judy{-- +(5mm,0mm) -- +(5mm,5mm) -- +(0mm,5mm) -- cycle [line width=0.6pt,fill=gainsboro,dotted]}
\tikzstyle{cc}=[circle,draw=black,fill=yellow, line width=0.5pt, inner sep=1.5]
\tikzstyle{rdot}=[circle,fill=red,draw=red,inner sep=1.5]

\draw (0mm,0mm)\hezi;
\draw (5mm,0mm)\hezi;
\draw (10mm,0mm)\hezi;
\draw (15mm,0mm)\judy;
\draw (20mm,0mm)\hezi;
\draw (0mm,-5mm)\hezi;
\draw (5mm,-5mm)\judy;
\draw (10mm,-5mm)\hezi;
\draw (15mm,-5mm)\hezi;
\draw (20mm,-5mm)\hezi;
\draw (0mm,-10mm)\hezi;
\draw (5mm,-10mm)\hezi;
\draw (10mm,-10mm)\judy;
\draw (15mm,-10mm)\hezi;
\draw (20mm,-10mm)\hezi;
\draw (0mm,-15mm)\judy;
\draw (5mm,-15mm)\hezi;
\draw (10mm,-15mm)\hezi;
\draw (15mm,-15mm)\hezi;
\draw (20mm,-15mm)\hezi;
\draw (0mm,-20mm)\hezi;
\draw (5mm,-20mm)\hezi;
\draw (10mm,-20mm)\hezi;
\draw (15mm,-20mm)\hezi;
\draw (20mm,-20mm)\judy;

\draw (0mm,5mm)--(15mm,5mm)--(20mm,0mm)--(25mm,0mm)[line width=1.25pt,draw=red];
\node[cc] at (0mm,5mm){};
\node at (-4mm,5mm){\red{$0$}};

\draw (0mm,0mm)--(5mm,0mm)--(15mm,-10mm)--(25mm,-10mm)[line width=1.25pt,draw=red];
\node[cc] at (15mm,-10mm){};
\node at (-4mm,0mm){\red{$0$}};

\draw (0mm,-10mm)--(5mm,-15mm)--(20mm,-15mm)--(25mm,-20mm)[line width=1.25pt,draw=red];
\node[cc] at (25mm,-20mm){};
\node at (-4mm,-10mm){\red{$0$}};


\end{tikzpicture}
\end{center}
By Theorem \ref{thmpaths},
we have
$\left|\mathcal{B}^{(1)}_{n-1,k}\right|=(k+1)I_{n-1,k}$ that is the \textbf{first} term on the right side of \eqref{recInk}.
For  $(\sigma,(i,i))\in\mathcal{B}^{(1)}_{n-1,k}$, since $(i,i)$ is of $d$-type $(0,0)$ by Proposition \ref{lemij=ji},
define
\begin{equation}\label{eqpfIn1st}
  \Psi((\sigma,(i,i)))=(\varphi_{(i,i)}(\sigma),i),
\end{equation}
and we have $(\varphi_{(i,i)}(\sigma),i)\in\widetilde{\mathcal{I}}_{n,k}$ by Lemma \ref{lemiip}.


Let
\begin{equation}\label{defB2nk}
  \begin{aligned}
\mathcal{B}^{(2)}_{n,k}=\{&(\sigma,(i,j))\mid\text{ $\sigma\in\mathcal{I}_{n,k}$,  and $(i,j)$ is the first grid point at which a $1_h$-path }\\
    &\quad\text{touches the main diagonal \emph{or} a filled square on the main diagonal of $P_\sigma$}\}.
  \end{aligned}
\end{equation}
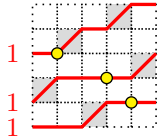
\begin{center}
\begin{tikzpicture}[scale = 0.65]
\def\hezi{-- +(5mm,0mm) -- +(5mm,5mm) -- +(0mm,5mm) -- cycle [line width=0.6pt,dotted]}
\def\judy{-- +(5mm,0mm) -- +(5mm,5mm) -- +(0mm,5mm) -- cycle [line width=0.6pt,fill=gainsboro,dotted]}
\tikzstyle{cc}=[circle,draw=black,fill=yellow, line width=0.5pt, inner sep=1.5]
\tikzstyle{rdot}=[circle,fill=red,draw=red,inner sep=1.5]

\draw (0mm,0mm)\hezi;
\draw (5mm,0mm)\hezi;
\draw (10mm,0mm)\hezi;
\draw (15mm,0mm)\judy;
\draw (20mm,0mm)\hezi;
\draw (0mm,-5mm)\hezi;
\draw (5mm,-5mm)\judy;
\draw (10mm,-5mm)\hezi;
\draw (15mm,-5mm)\hezi;
\draw (20mm,-5mm)\hezi;
\draw (0mm,-10mm)\hezi;
\draw (5mm,-10mm)\hezi;
\draw (10mm,-10mm)\hezi;
\draw (15mm,-10mm)\hezi;
\draw (20mm,-10mm)\judy;
\draw (0mm,-15mm)\judy;
\draw (5mm,-15mm)\hezi;
\draw (10mm,-15mm)\hezi;
\draw (15mm,-15mm)\hezi;
\draw (20mm,-15mm)\hezi;
\draw (0mm,-20mm)\hezi;
\draw (5mm,-20mm)\hezi;
\draw (10mm,-20mm)\judy;
\draw (15mm,-20mm)\hezi;
\draw (20mm,-20mm)\hezi;

\draw (0mm,-5mm)--(5mm,-5mm)--(10mm,0mm)--(15mm,0mm)--(20mm,5mm)--(25mm,5mm)[line width=1.25pt,draw=red];
\node[cc] at (5mm,-5mm){};
\node at (-4mm,-5mm){\red{$1$}};

\draw (0mm,-15mm)--(5mm,-10mm)--(20mm,-10mm)--(25mm,-5mm)[line width=1.25pt,draw=red];
\node[cc] at (15mm,-10mm){};
\node at (-4mm,-15mm){\red{$1$}};

\draw (0mm,-20mm)--(10mm,-20mm)--(15mm,-15mm)--(25mm,-15mm)[line width=1.25pt,draw=red];
\node[cc] at (20mm,-15mm){};
\node at (-4mm,-20mm){\red{$1$}};

\end{tikzpicture}
\end{center}
It follows from Theorem \ref{thmpaths} that $\left|\mathcal{B}^{(2)}_{n-1,k-1}\right|=(n-k)I_{n-1,k-1}$,
which is the \textbf{second} term on the right side of \eqref{recInk}.
For $(\sigma,(i,j))\in \mathcal{B}^{(2)}_{n-1,k-1}$,
set
\begin{equation}\label{eqpfIn2nd}
\Psi((\sigma, (i,j)))=(\varphi_{(i,j)}(\sigma),i).
\end{equation}
Since $i=j$ \emph{or} $i=j+1$ with $\langle j,j\rangle$ filled in $P_\sigma$,
we deduce $(\varphi_{(i,j)}(\sigma),i)\in\widetilde{\mathcal{I}}_{n,k}$ by Propositions \ref{lemiip} and \ref{lemi+1i}.
%
%

Let
\begin{equation}\label{defB3nk}
\mathcal{B}^{(3)}_{n,k}=\{(\sigma,(i,j))\mid \sigma\in \mathcal{I}_{n,k}\text{ and $(i,j)$ is of $d$-type $(0,0)$ in $P_\sigma$}\},
\end{equation}
From Theorem \ref{thmdtype},
We see $\left|\mathcal{B}^{(3)}_{n-2,k}\right|=\left((k+1)^2+n-2\right)I_{n-2,k}$,
which is the \textbf{third} term on the right side of \eqref{recInk}.
For $(\sigma,(i,j))\in \mathcal{B}^{(3)}_{n-2,k}$,
define
\begin{equation}\label{eqpfIn3rd}
\Psi((\sigma,(i,j)))=\left\{
        \begin{array}{ll}
          \left(\xi_{(i,j)}(\sigma),\chi_{ij}\right), & i\neq j; \\[6pt]
          \left(\eta_i(\sigma),i\right), & i=j,
        \end{array}
      \right.
\end{equation}
where
\begin{equation*}
\chi_{ij}=\left\{
\begin{array}{ll}
i, & \hbox{if  $i< j$,} \\[6pt]
i+1, & \hbox{if $i>j$.}
\end{array}
    \right.
\end{equation*}
By setting $(p,q)=(0,0)$ in Proposition \ref{lemxi} and $p=0$ in Proposition \ref{lemetaeta'} for the operation $\eta_i$,
we obtain $\Psi(\sigma,(i,j))\in\widetilde{\mathcal{I}}_{n,k}$.

Let
\begin{equation}\label{defB4nk}
  \begin{aligned}
\mathcal{B}^{(4)}_{n,k}=\{(\sigma,(i,j))\mid &\text{ $\sigma\in\mathcal{I}_{n,k}$,  and $(i,j)$ is the grid point }\\
&\qquad\text{with $i=j$ \emph{or} of $d$-type $(0,1)$ or $(1,0)$ in $P_\sigma$}\}.
  \end{aligned}
\end{equation}
For $\sigma\in \mathcal{I}_{n-2,k-1}$,
since the grid points of $d$-type $(0,1)$ and $(1,0)$
are not on the main diagonal by Proposition \ref{lemij=ji},
with Theorem \ref{thmdtype}, we get
\begin{align*}
\left|\mathcal{B}^{(4)}_{n-2,k-1}\right|&=
(n-1 + 2(k(n-k-1)-(n-2)))I_{n-2,k-1}\\[4pt]
&=
(2 k(n-k-1)-n+3)I_{n-2,k-1},
\end{align*}
which is the \textbf{fourth} term on the right side of \eqref{recInk}.
For $(\sigma,(i,j))\in \mathcal{B}^{(4)}_{n-2,k-1}$,
let
\begin{equation}\label{eqpfIn4th}
\Psi((\sigma,(i,j)))=\left\{
        \begin{array}{ll}
          \left(\xi_{(i,j)}(\sigma),\chi_{ij}\right), & i\neq j; \\[6pt]
          \left(\eta_i(\sigma),i\right), & i=j \text{ and $(i,j)$ is of $d$-type $(1,1)$}; \\[6pt]
          \left(\eta'_i(\sigma),i\right), & i=j \text{ and $(i,j)$ is of $d$-type $(0,0)$}.
        \end{array}
      \right.
\end{equation}
Note that by Proposition \ref{lemij=ji},
the grid point $(i,j)$ with $i=j$ must have $d$-type $(0,0)$ or $(1,1)$,
and the grid points $(i,j)$ with $i\neq j$ must have $d$-type $(0,1)$ or $(1,0)$.
Thus from Propositions \ref{lemxi} and \ref{lemetaeta'},
we see 
$\Psi((\sigma,(i,j)))\in\widetilde{\mathcal{I}}_{n,k}$.

Let
\begin{equation}\label{defB5nk}
\mathcal{B}^{(5)}_{n,k}=\{(\sigma,(i,j))\mid \sigma\in \mathcal{I}_{n,k}\text{ and $(i,j)$ is of $d$-type $(1,1)$ in $P_\sigma$}\}.
\end{equation}
By Theorem \ref{thmdtype}, we deduce
$\left|\mathcal{B}^{(5)}_{n-2,k-2}\right|=\left((n-k)^{2}+n-2\right) I_{n-2, k-2}$,
which is the \textbf{fifth} term on the right side of \eqref{recInk}.
For  $(\sigma,(i,j))\in \mathcal{B}^{(5)}_{n-2,k-2}$,
set
\begin{equation}\label{eqpfIn5th}
\Psi((\sigma,(i,j)))=\left\{
        \begin{array}{ll}
          \left(\xi_{(i,j)}(\sigma),\chi_{ij}\right), & i\neq j; \\[6pt]
          \left(\eta'_i(\sigma),i\right), & i=j.
        \end{array}
      \right.
\end{equation}
Using Proposition \ref{lemxi} with $(p,q)=(1,1)$ and Proposition $\ref{lemetaeta'}$ with $p=1$ for the operation $\eta_i'$,
we see 
$\Psi(\sigma,(i,j))\in\widetilde{\mathcal{I}}_{n,k}$.

It is clear that the above five sets $\mathcal{B}^{(1)}_{n-1,k}$,
$\mathcal{B}^{(2)}_{n-1,k-1}$, $\mathcal{B}^{(3)}_{n-2,k}$, $\mathcal{B}^{(4)}_{n-2,k-1}$, $\mathcal{B}^{(5)}_{n-2,k-2}$ are pairwise disjoint,
which implies that the cardinality of the set
\begin{equation}\label{defBnk}
\widetilde{\mathcal{B}}_{n,k}:=
\mathcal{B}^{(1)}_{n-1,k}\uplus\mathcal{B}^{(2)}_{n-1,k-1}\uplus\mathcal{B}^{(3)}_{n-2,k}
\uplus\mathcal{B}^{(4)}_{n-2,k-1}\uplus\mathcal{B}^{(5)}_{n-2,k-2}.
\end{equation}
exactly equals the right side of \eqref{recInk}.
Moreover, combining  \eqref{eqpfIn1st}, \eqref{eqpfIn2nd}, \eqref{eqpfIn3rd}, \eqref{eqpfIn4th} and \eqref{eqpfIn5th},
we establish the mapping
\[
\Psi\colon \widetilde{\mathcal{B}}_{n,k}\rightarrow \widetilde{\mathcal{I}}_{n,k}.
\]

To complete our proof,
for $(\pi,i)\in\widetilde{\mathcal{I}}_{n,k}$, let
\begin{equation}\label{theta}
\Theta((\pi,i))=\left\{
        \begin{array}{ll}
          \left(\xi^{-1}_{\langle i,\pi_i-1\rangle}(\pi),\left(i,\pi_i-1\right)\right), & \pi_i>i+1; \\[8pt]
          \left(\xi^{-1}_{\langle i-1,\pi_i\rangle}(\pi),\left(i-1,\pi_i\right)\right), & \pi_i<i-1; \\[8pt]
          \left(\eta'^{-1}_i(\pi),\left(i,\pi_i-1\right)\right), & \pi_i=i+1; \\[8pt]
          \left(\varphi^{-1}_{\langle i,\pi_i\rangle}(\pi),\left(i,\pi_i\right)\right), & \pi_i=i-1; \\[8pt]
          \left(\varphi^{-1}_{\langle i,\pi_i\rangle}(\pi),\left(i,\pi_i\right)\right), & \pi_i=i\text{ and $\pi_{i+1}\neq i+1$}; \\[8pt]
          \left(\eta^{-1}_{i}(\pi),\left(i,\pi_{i}\right)\right), & \pi_i=i\text{ and $\pi_{i+1}=i+1$}.
        \end{array}
      \right.
\end{equation}
We shall show that $\Theta$ is the inverse of $\Psi$.

\noindent{\bf{Case 1: $\bm{\pi_i>i+1}$.}}
By the definition of $\xi^{-1}_{\ij}$ with $i<j$
and Proposition \ref{lemxi},
we see
$$\xi^{-1}_{{\langle i,\pi_i-1\rangle}}(\pi)\in\mathcal{I}_{n-2,k-p-q}$$ if $(i,\pi_i-1)$ is of $d$-type $(p,q)$ in the grid of $\xi^{-1}_{{\langle i,\pi_i-1\rangle}}(\pi)$,
where $p,q\in\{0,1\}$.
See the following grids for an example.
\begin{center}
\begin{tikzpicture}[scale = 0.65]
\def\hezi{-- +(5mm,0mm) -- +(5mm,5mm) -- +(0mm,5mm) -- cycle [line width=0.6pt]}
\def\judy{-- +(5mm,0mm) -- +(5mm,5mm) -- +(0mm,5mm) -- cycle [line width=0.6pt,fill=gainsboro]}
\tikzstyle{cc}=[circle,draw=black,fill=yellow, line width=0.5pt, inner sep=1.5]
\tikzstyle{dd}=[circle,double=yellow,draw=black,fill=white, line width=0.5pt, inner sep=1.4]
\tikzstyle{rdot}=[circle,fill=red,draw=red,inner sep=2]
\tikzstyle{bdot}=[circle,fill=babyblue,draw=babyblue,inner sep=2]

\draw (0mm,0mm)\hezi;
\draw (5mm,0mm)\hezi;
\draw (10mm,0mm)\hezi;
\draw (15mm,0mm)-- +(5mm,0mm) -- +(5mm,5mm) -- +(0mm,5mm) -- cycle [line width=0.8pt,fill=red];
\draw (20mm,0mm)\hezi;
\draw (0mm,-5mm)\hezi;
\draw (5mm,-5mm)\judy;
\draw (10mm,-5mm)\hezi;
\draw (15mm,-5mm)\hezi;
\draw (20mm,-5mm)\hezi;
\draw (0mm,-10mm)\hezi;
\draw (5mm,-10mm)\hezi;
\draw (10mm,-10mm)\hezi;
\draw (15mm,-10mm)\hezi;
\draw (20mm,-10mm)\judy;
\draw (0mm,-15mm)\judy;
\draw (5mm,-15mm)\hezi;
\draw (10mm,-15mm)\hezi;
\draw (15mm,-15mm)\hezi;
\draw (20mm,-15mm)\hezi;
\draw (0mm,-20mm)\hezi;
\draw (5mm,-20mm)\hezi;
\draw (10mm,-20mm)\judy;
\draw (15mm,-20mm)\hezi;
\draw (20mm,-20mm)\hezi;



\node at (12.5mm,-26mm){$42513\in\mathcal{I}_{5,2}$};


\begin{scope}[shift={(60mm,-6.5mm)},scale=1]
\draw (0mm,0mm)\judy;
\draw (5mm,0mm)\hezi;
\draw (10mm,0mm)\hezi;
\draw (0mm,-5mm)\hezi;
\draw (5mm,-5mm)\hezi;
\draw (10mm,-5mm)\judy;
\draw (0mm,-10mm)\hezi;
\draw (5mm,-10mm)\judy;
\draw (10mm,-10mm)\hezi;

\node[rdot]at(10mm,5mm){};

\node at (7.5mm,-16mm){$132\in\mathcal{I}_{3,1}$};

\end{scope}
\end{tikzpicture}
\end{center}
Thus by setting
\begin{align*}
  &\mathcal{B}^{(3)_1}_{n,k}=\{(\sigma,(i,j))\mid \sigma\in \mathcal{I}_{n,k}\text{ and $(i,j)$ is of $d$-type $(0,0)$ in $P_\sigma$ with $i<j$}\},\\
  &\mathcal{B}^{(4)_1}_{n,k}=\{(\sigma,(i,j))\mid \sigma\in \mathcal{I}_{n,k}\text{ and  $(i,j)$ is of $d$-type $(0,1)$ or $(1,0)$ in $P_\sigma$ with $i<j$}\}, \\
  &\mathcal{B}^{(5)_1}_{n,k}=\{(\sigma,(i,j))\mid \sigma\in \mathcal{I}_{n,k}\text{ and  $(i,j)$ is of $d$-type $(1,1)$ in $P_\sigma$ with $i<j$}\},
\end{align*}
we deduce from \eqref{defB3nk}, \eqref{defB4nk} and \eqref{defB5nk} that $\mathcal{B}^{(r)_1}_{n,k}\subseteq \mathcal{B}^{(r)}_{n,k}$ for $r\in\{3,4,5\}$, and
\begin{equation*}
  \left(\xi^{-1}_{\langle i,\pi_i-1\rangle}(\pi),\left(i,\pi_i-1\right)\right)\in \mathcal{B}^{(3)_1}_{n-2,k}\uplus\mathcal{B}^{(4)_1}_{n-2,k-1}
  \uplus\mathcal{B}^{(5)_1}_{n-2,k-2}
  \subseteq \widetilde{\mathcal{B}}_{n,k}.
\end{equation*}

\noindent{\bf{Case 2: $\bm{\pi_i<i-1}$.}}
By the definition of $\xi^{-1}_{\ij}$ with $i>j$
and Proposition \ref{lemxi},
if the grid point $(i-1,\pi_i)$ is of $d$-type $(p,q)$ for $p,q\in\{0,1\}$
in the grid of $\xi^{-1}_{\langle i-1,\pi_i\rangle}(\pi)$, then
we have
$$\xi^{-1}_{\langle i-1,\pi_i\rangle}(\pi)\in\mathcal{I}_{n-2,k-p-q}.$$
Refer the following grids  for a clearer view.
\begin{center}
\begin{tikzpicture}[scale = 0.65]
\def\hezi{-- +(5mm,0mm) -- +(5mm,5mm) -- +(0mm,5mm) -- cycle [line width=0.6pt]}
\def\judy{-- +(5mm,0mm) -- +(5mm,5mm) -- +(0mm,5mm) -- cycle [line width=0.6pt,fill=gainsboro]}
\tikzstyle{cc}=[circle,draw=black,fill=yellow, line width=0.5pt, inner sep=1.5]
\tikzstyle{dd}=[circle,double=yellow,draw=black,fill=white, line width=0.5pt, inner sep=1.4]
\tikzstyle{rdot}=[circle,fill=red,draw=red,inner sep=2]
\tikzstyle{bdot}=[circle,fill=babyblue,draw=babyblue,inner sep=2]

\draw (0mm,0mm)\hezi;
\draw (5mm,0mm)\hezi;
\draw (10mm,0mm)\hezi;
\draw (15mm,0mm)\judy;
\draw (20mm,0mm)\hezi;
\draw (0mm,-5mm)\hezi;
\draw (5mm,-5mm)\judy;
\draw (10mm,-5mm)\hezi;
\draw (15mm,-5mm)\hezi;
\draw (20mm,-5mm)\hezi;
\draw (0mm,-10mm)\hezi;
\draw (5mm,-10mm)\hezi;
\draw (10mm,-10mm)\hezi;
\draw (15mm,-10mm)\hezi;
\draw (20mm,-10mm)\judy;
\draw (0mm,-15mm)\judy;
\draw (5mm,-15mm)\hezi;
\draw (10mm,-15mm)\hezi;
\draw (15mm,-15mm)\hezi;
\draw (20mm,-15mm)\hezi;
\draw (0mm,-20mm)\hezi;
\draw (5mm,-20mm)\hezi;
\draw (10mm,-20mm)-- +(5mm,0mm) -- +(5mm,5mm) -- +(0mm,5mm) -- cycle [line width=0.8pt,fill=red];
\draw (15mm,-20mm)\hezi;
\draw (20mm,-20mm)\hezi;



\node at (12.5mm,-26mm){$42513\in\mathcal{I}_{5,2}$};


\begin{scope}[shift={(60mm,-6.5mm)},scale=1]
\draw (0mm,0mm)\hezi;
\draw (5mm,0mm)\hezi;
\draw (10mm,0mm)\judy;
\draw (0mm,-5mm)\hezi;
\draw (5mm,-5mm)\judy;
\draw (10mm,-5mm)\hezi;
\draw (0mm,-10mm)\judy;
\draw (5mm,-10mm)\hezi;
\draw (10mm,-10mm)\hezi;

\node[rdot]at(10mm,-10mm){};

\node at (7.5mm,-16mm){$321\in\mathcal{I}_{3,2}$};

\end{scope}


%
%

\end{tikzpicture}
\end{center}
Hence by letting
\begin{align*}
  &\mathcal{B}^{(3)_2}_{n,k}=\{(\sigma,(i,j))\mid \sigma\in \mathcal{I}_{n,k}\text{ and  $(i,j)$ is of $d$-type $(0,0)$ in $P_\sigma$ with $i>j$}\},\\
  &\mathcal{B}^{(4)_2}_{n,k}=\{(\sigma,(i,j))\mid \sigma\in \mathcal{I}_{n,k}\text{ and $(i,j)$ is of $d$-type $(0,1)$ or $(1,0)$ in $P_\sigma$ with $i>j$}\}, \\
  &\mathcal{B}^{(5)_2}_{n,k}=\{(\sigma,(i,j))\mid \sigma\in \mathcal{I}_{n,k}\text{ and  $(i,j)$ is of $d$-type $(1,1)$ in $P_\sigma$ with $i>j$}\},
\end{align*}
we have $\mathcal{B}^{(r)_2}_{n,k}\subseteq \mathcal{B}^{(r)}_{n,k}$ for $r\in\{3,4,5\}$ by \eqref{defB3nk}, \eqref{defB4nk} and \eqref{defB5nk},
which leads to
\begin{equation*}
  \left(\xi^{-1}_{\langle i-1,\pi_i\rangle}(\pi),\left(i-1,\pi_i\right)\right)\in \mathcal{B}^{(3)_2}_{n-2,k}\uplus\mathcal{B}^{(4)_2}_{n-2,k-1}
  \uplus\mathcal{B}^{(5)_2}_{n-2,k-2}
  \subseteq \widetilde{\mathcal{B}}_{n,k}.
\end{equation*}

\noindent{\bf{Case 3: $\bm{\pi_i=i+1}$.}}
This
implies that the squares $\langle i,i+1 \rangle$ and $\langle i+1,i\rangle$ are both filled in $P_\pi$.
Thus by Proposition \ref{lemetaeta'},
we see
$$\eta'^{-1}_i(\pi)\in \mathcal{I}_{n-2,k-1-p}$$
if the grid point $(i,i)$ is of $d$-type $(p,p)$ in the grid of $\eta'^{-1}_i(\pi)$,
as demonstrated by the next example.
\begin{center}
\begin{tikzpicture}[scale = 0.65]
\def\hezi{-- +(5mm,0mm) -- +(5mm,5mm) -- +(0mm,5mm) -- cycle [line width=0.6pt]}
\def\judy{-- +(5mm,0mm) -- +(5mm,5mm) -- +(0mm,5mm) -- cycle [line width=0.6pt,fill=gainsboro]}
\tikzstyle{cc}=[circle,draw=black,fill=yellow, line width=0.5pt, inner sep=1.5]
\tikzstyle{dd}=[circle,double=yellow,draw=black,fill=white, line width=0.5pt, inner sep=1.4]
\tikzstyle{rdot}=[circle,fill=red,draw=red,inner sep=2]
\tikzstyle{bdot}=[circle,fill=babyblue,draw=babyblue,inner sep=2]

\draw (0mm,0mm)\judy;
\draw (5mm,0mm)\hezi;
\draw (10mm,0mm)\hezi;
\draw (15mm,0mm)\hezi;
\draw (20mm,0mm)\hezi;
\draw (0mm,-5mm)\hezi;
\draw (5mm,-5mm)\hezi;
\draw (10mm,-5mm)-- +(5mm,0mm) -- +(5mm,5mm) -- +(0mm,5mm) -- cycle [line width=0.8pt,fill=red];
\draw (15mm,-5mm)\hezi;
\draw (20mm,-5mm)\hezi;
\draw (0mm,-10mm)\hezi;
\draw (5mm,-10mm)\judy;
\draw (10mm,-10mm)\hezi;
\draw (15mm,-10mm)\hezi;
\draw (20mm,-10mm)\hezi;
\draw (0mm,-15mm)\hezi;
\draw (5mm,-15mm)\hezi;
\draw (10mm,-15mm)\hezi;
\draw (15mm,-15mm)\hezi;
\draw (20mm,-15mm)\judy;
\draw (0mm,-20mm)\hezi;
\draw (5mm,-20mm)\hezi;
\draw (10mm,-20mm)\hezi;
\draw (15mm,-20mm)\judy;
\draw (20mm,-20mm)\hezi;
\node at (12.5mm,-26mm){$13254\in\mathcal{I}_{5,2}$};


\begin{scope}[shift={(60mm,-6.5mm)},scale=1]
\draw (0mm,0mm)\judy;
\draw (5mm,0mm)\hezi;
\draw (10mm,0mm)\hezi;
\draw (0mm,-5mm)\hezi;
\draw (5mm,-5mm)\hezi;
\draw (10mm,-5mm)\judy;
\draw (0mm,-10mm)\hezi;
\draw (5mm,-10mm)\judy;
\draw (10mm,-10mm)\hezi;

\node[rdot] at(5mm,0mm){};

\node at (7.5mm,-16mm){$132\in\mathcal{I}_{3,1}$};

\end{scope}


\end{tikzpicture}
\end{center}
Let
\begin{equation}\label{defB43nk}
  \mathcal{B}^{(4)_3}_{n,k}=\{(\sigma,(i,i))\mid \sigma\in \mathcal{I}_{n,k}\text{ and  $(i,i)$ is of $d$-type $(0,0)$ in $P_\sigma$}\}
\end{equation}
and
\begin{equation}\label{defB53nk}
  \mathcal{B}^{(5)_3}_{n,k}=\{(\sigma,(i,i))\mid \sigma\in \mathcal{I}_{n,k}\text{ and  $(i,i)$ is of $d$-type $(1,1)$ in $P_\sigma$}\},
\end{equation}
then $\mathcal{B}^{(r)_3}_{n,k}\subseteq \mathcal{B}^{(r)}_{n,k}$ for $r\in\{4,5\}$
by \eqref{defB4nk} and \eqref{defB5nk}.
Hence, we have
\begin{equation*}
  \left(\eta'^{-1}_i(\pi),\left(i,\pi_i-1\right)\right)\in \mathcal{B}^{(4)_3}_{n-2,k-1}\uplus\mathcal{B}^{(5)_3}_{n-2,k-2}
  \subseteq \widetilde{\mathcal{B}}_{n,k}.
\end{equation*}

\noindent{\bf{Case 4: $\bm {\pi_i=i-1}$.}}
Similar to Case 3,
the squares
$\langle i,i-1\rangle$ and $\langle i-1,i\rangle$ are both filled in  $P_\pi$.
Thus by  Proposition \ref{lemi+1i},
deleting the square $\langle i,\pi_i \rangle$ from $P_\pi$ yields
$\varphi^{-1}_{\langle i,\pi_i\rangle}(\pi)\in\mathcal{I}_{n-1}$,
and we have
$$\varphi^{-1}_{\langle i,\pi_i\rangle}(\pi)\in\mathcal{I}_{n-1,k-1}.$$
In the following example,
one should note that
$(i,\pi_i)$ must be the grid point at which a certain $1_h$-path first touches
the filled square $\langle i-1,\pi_i\rangle$ on the main diagonal
of the grid of  $\varphi^{-1}_{\langle i,\pi_i\rangle}(\pi)$.
\begin{center}
\begin{tikzpicture}[scale = 0.65]
\def\hezi{-- +(5mm,0mm) -- +(5mm,5mm) -- +(0mm,5mm) -- cycle [line width=0.6pt]}
\def\judy{-- +(5mm,0mm) -- +(5mm,5mm) -- +(0mm,5mm) -- cycle [line width=0.6pt,fill=gainsboro]}
\tikzstyle{cc}=[circle,draw=black,fill=yellow, line width=0.5pt, inner sep=1.5]
\tikzstyle{dd}=[circle,double=yellow,draw=black,fill=white, line width=0.5pt, inner sep=1.4]
\tikzstyle{rdot}=[circle,fill=red,draw=red,inner sep=2]
\tikzstyle{bdot}=[circle,fill=babyblue,draw=babyblue,inner sep=2]

\draw (0mm,0mm)\judy;
\draw (5mm,0mm)\hezi;
\draw (10mm,0mm)\hezi;
\draw (15mm,0mm)\hezi;
\draw (20mm,0mm)\hezi;
\draw (0mm,-5mm)\hezi;
\draw (5mm,-5mm)\hezi;
\draw (10mm,-5mm)\judy;
\draw (15mm,-5mm)\hezi;
\draw (20mm,-5mm)\hezi;
\draw (0mm,-10mm)\hezi;
\draw (5mm,-10mm)-- +(5mm,0mm) -- +(5mm,5mm) -- +(0mm,5mm) -- cycle [line width=0.8pt,fill=red];
\draw (10mm,-10mm)\hezi;
\draw (15mm,-10mm)\hezi;
\draw (20mm,-10mm)\hezi;
\draw (0mm,-15mm)\hezi;
\draw (5mm,-15mm)\hezi;
\draw (10mm,-15mm)\hezi;
\draw (15mm,-15mm)\hezi;
\draw (20mm,-15mm)\judy;
\draw (0mm,-20mm)\hezi;
\draw (5mm,-20mm)\hezi;
\draw (10mm,-20mm)\hezi;
\draw (15mm,-20mm)\judy;
\draw (20mm,-20mm)\hezi;

\node at (12.5mm,-26mm){$13254\in\mathcal{I}_{5,2}$};


\begin{scope}[shift={(60mm,-5mm)},scale=1]
\draw (0mm,0mm)\judy;
\draw (5mm,0mm)\hezi;
\draw (10mm,0mm)\hezi;
\draw (15mm,0mm)\hezi;
\draw (0mm,-5mm)\hezi;
\draw (5mm,-5mm)\judy;
\draw (10mm,-5mm)\hezi;
\draw (15mm,-5mm)\hezi;
\draw (0mm,-10mm)\hezi;
\draw (5mm,-10mm)\hezi;
\draw (10mm,-10mm)\hezi;
\draw (15mm,-10mm)\judy;
\draw (0mm,-15mm)\hezi;
\draw (5mm,-15mm)\hezi;
\draw (10mm,-15mm)\judy;
\draw (15mm,-15mm)\hezi;

\node[rdot] at(5mm,-5mm){};

\node at (10mm,-21mm){$1243\in\mathcal{I}_{4,1}$};

\end{scope}


\end{tikzpicture}
\end{center}
Denote by
\begin{equation*}
  \begin{aligned}
\mathcal{B}^{(2)_1}_{n,k}=\{(\sigma,(i,j))\mid&\text{ $\sigma\in\mathcal{I}_{n,k}$,  and $(i,j)$ is the first grid point at which }\\
    &\quad\text{a $1_h$-path touches a filled square on the main diagonal of $P_\sigma$}\}.
  \end{aligned}
\end{equation*}
It follows from \eqref{defB2nk} that $\mathcal{B}^{(2)_1}_{n,k}\subseteq \mathcal{B}^{(2)}_{n,k}$,
and
\begin{equation*}
  \left(\varphi^{-1}_{\langle i,\pi_i\rangle}(\pi),\left(i,\pi_i\right)\right)\in \mathcal{B}^{(2_1)}_{n-1,k-1}
  \subseteq \widetilde{\mathcal{B}}_{n,k}.
\end{equation*}

\noindent{\bf{Case 5: $\bm {\pi_i=i}$ and $\bm{\pi_{i+1}\neq i+1}$.}}
By the fact $\pi=i$, we have
\[
\varphi^{-1}_{\langle i,\pi_i\rangle}(\pi)\in\mathcal{I}_{n-1,k-p},
\]
if $(i,\pi_i)$ is of $d$-type $(p,p)$ for $p\in\{0,1\}$
in the grid of $\varphi^{-1}_{\langle i,\pi_i\rangle}(\pi)$.
Since the square $\langle i+1,\pi_i+1\rangle$ is not filled in $P_\pi$,
the grid point $(i,\pi_i)$ in the grid of $\varphi^{-1}_{\langle i,\pi_i\rangle}(\pi)$ must be the point at which
a $0_h$-path last touches the main diagonal,
or the only point at which a $1_h$-path touches the main diagonal.
See the following figures for an example.

\begin{center}
\begin{tikzpicture}[scale = 0.65]
\def\hezi{-- +(5mm,0mm) -- +(5mm,5mm) -- +(0mm,5mm) -- cycle [line width=0.6pt]}
\def\judy{-- +(5mm,0mm) -- +(5mm,5mm) -- +(0mm,5mm) -- cycle [line width=0.6pt,fill=gainsboro]}
\tikzstyle{cc}=[circle,draw=black,fill=yellow, line width=0.5pt, inner sep=1.5]
\tikzstyle{dd}=[circle,double=yellow,draw=black,fill=white, line width=0.5pt, inner sep=1.4]
\tikzstyle{rdot}=[circle,fill=red,draw=red,inner sep=2]
\tikzstyle{bdot}=[circle,fill=babyblue,draw=babyblue,inner sep=2]

\draw (0mm,0mm)\judy;
\draw (5mm,0mm)\hezi;
\draw (10mm,0mm)\hezi;
\draw (15mm,0mm)\hezi;
\draw (20mm,0mm)\hezi;
\draw (0mm,-5mm)\hezi;
\draw (5mm,-5mm)-- +(5mm,0mm) -- +(5mm,5mm) -- +(0mm,5mm) -- cycle [line width=0.8pt,fill=babyblue];
\draw (10mm,-5mm)\hezi;
\draw (15mm,-5mm)\hezi;
\draw (20mm,-5mm)\hezi;
\draw (0mm,-10mm)\hezi;
\draw (5mm,-10mm)\hezi;
\draw (10mm,-10mm)\hezi;
\draw (15mm,-10mm)\hezi;
\draw (20mm,-10mm)\judy;
\draw (0mm,-15mm)\hezi;
\draw (5mm,-15mm)\hezi;
\draw (10mm,-15mm)\hezi;
\draw (15mm,-15mm)-- +(5mm,0mm) -- +(5mm,5mm) -- +(0mm,5mm) -- cycle [line width=0.8pt,fill=red];
\draw (20mm,-15mm)\hezi;
\draw (0mm,-20mm)\hezi;
\draw (5mm,-20mm)\hezi;
\draw (10mm,-20mm)\judy;
\draw (15mm,-20mm)\hezi;
\draw (20mm,-20mm)\hezi;

\node at (12.5mm,-26mm){$12543\in\mathcal{I}_{5,2}$};


\begin{scope}[shift={(50mm,-5mm)}]
\draw (0mm,0mm)\judy;
\draw (5mm,0mm)\hezi;
\draw (10mm,0mm)\hezi;
\draw (15mm,0mm)\hezi;
\draw (0mm,-5mm)\hezi;
\draw (5mm,-5mm)\hezi;
\draw (10mm,-5mm)\hezi;
\draw (15mm,-5mm)\judy;
\draw (0mm,-10mm)\hezi;
\draw (5mm,-10mm)\hezi;
\draw (10mm,-10mm)\judy;
\draw (15mm,-10mm)\hezi;
\draw (0mm,-15mm)\hezi;
\draw (5mm,-15mm)\judy;
\draw (10mm,-15mm)\hezi;
\draw (15mm,-15mm)\hezi;

\node[bdot] at(5mm,0mm){};

\node at (10mm,-21mm){$1432\in\mathcal{I}_{4,2}$};

\end{scope}


\begin{scope}[shift={(100mm,-5mm)}]
\draw (0mm,0mm)\judy;
\draw (5mm,0mm)\hezi;
\draw (10mm,0mm)\hezi;
\draw (15mm,0mm)\hezi;
\draw (0mm,-5mm)\hezi;
\draw (5mm,-5mm)\judy;
\draw (10mm,-5mm)\hezi;
\draw (15mm,-5mm)\hezi;
\draw (0mm,-10mm)\hezi;
\draw (5mm,-10mm)\hezi;
\draw (10mm,-10mm)\hezi;
\draw (15mm,-10mm)\judy;
\draw (0mm,-15mm)\hezi;
\draw (5mm,-15mm)\hezi;
\draw (10mm,-15mm)\judy;
\draw (15mm,-15mm)\hezi;

\node[rdot] at(15mm,-10mm){};

\node at (10mm,-21mm){$1243\in\mathcal{I}_{4,1}$};

\end{scope}

\end{tikzpicture}
\end{center}
Let
\begin{equation*}
  \begin{aligned}
\mathcal{B}^{(2)_2}_{n,k}=\{(\sigma,(i,j))\mid&\text{ $\sigma\in\mathcal{I}_{n,k}$,  and $(i,j)$ is the first grid point}\\
    &\quad\text{at which a $1_h$-path touches the main diagonal of $P_\sigma$}\},
  \end{aligned}
\end{equation*}
then we have $\mathcal{B}^{(2)_2}_{n,k}\subseteq \mathcal{B}^{(2)}_{n,k}$ by \eqref{defB2nk}.
Combing the definition \eqref{defB1nk} of the set $\mathcal{B}^{(1)}_{n,k}$,
we obtain that
\begin{equation*}
  \left(\varphi^{-1}_{\langle i,\pi_i\rangle}(\pi),\left(i,\pi_i\right)\right)\in \mathcal{B}^{(1)}_{n-1,k}\uplus \mathcal{B}^{(2_2)}_{n-1,k-1}
  \subseteq \widetilde{\mathcal{B}}_{n,k}.
\end{equation*}

\noindent{\bf{Case 6: $\bm {\pi_i=i}$ and $\bm{\pi_{i+1}=i+1}$.}}
Note that the  squares $\langle i,i\rangle$ and $\langle i+1,i+1\rangle$
are both filled in $P_\pi$,
which allows us to employ the double deleting operation $\eta^{-1}_{i}$.
Due to Proposition \ref{lemetaeta'}, we conclude that
\begin{equation*}
  \eta^{-1}_{i}(\pi)\in\mathcal{I}_{n-2,k-p}
\end{equation*}
if the grid point $(i,\pi_i)$ is of $d$-type $(p,p)$ in
the grid of $\eta^{-1}_{i}(\pi)$ for $p\in\{0,1\}$.
The following grids give an example.
\begin{center}
\begin{tikzpicture}[scale = 0.65]
\def\hezi{-- +(5mm,0mm) -- +(5mm,5mm) -- +(0mm,5mm) -- cycle [line width=0.6pt]}
\def\judy{-- +(5mm,0mm) -- +(5mm,5mm) -- +(0mm,5mm) -- cycle [line width=0.6pt,fill=gainsboro]}
\tikzstyle{cc}=[circle,draw=black,fill=yellow, line width=0.5pt, inner sep=1.5]
\tikzstyle{dd}=[circle,double=yellow,draw=black,fill=white, line width=0.5pt, inner sep=1.4]
\tikzstyle{rdot}=[circle,fill=red,draw=red,inner sep=2]
\tikzstyle{bdot}=[circle,fill=babyblue,draw=babyblue,inner sep=2]

\draw (0mm,0mm)\hezi;
\draw (5mm,0mm)\hezi;
\draw (10mm,0mm)\hezi;
\draw (15mm,0mm)\judy;
\draw (20mm,0mm)\hezi;
\draw (0mm,-5mm)\hezi;
\draw (5mm,-5mm)-- +(5mm,0mm) -- +(5mm,5mm) -- +(0mm,5mm) -- cycle [line width=0.8pt,fill=red];
\draw (10mm,-5mm)\hezi;
\draw (15mm,-5mm)\hezi;
\draw (20mm,-5mm)\hezi;
\draw (0mm,-10mm)\hezi;
\draw (5mm,-10mm)\hezi;
\draw (10mm,-10mm)\judy;
\draw (15mm,-10mm)\hezi;
\draw (20mm,-10mm)\hezi;
\draw (0mm,-15mm)\judy;
\draw (5mm,-15mm)\hezi;
\draw (10mm,-15mm)\hezi;
\draw (15mm,-15mm)\hezi;
\draw (20mm,-15mm)\hezi;
\draw (0mm,-20mm)\hezi;
\draw (5mm,-20mm)\hezi;
\draw (10mm,-20mm)\hezi;
\draw (15mm,-20mm)\hezi;
\draw (20mm,-20mm)\judy;


\node at (12.5mm,-26mm){$42315\in\mathcal{I}_{5,2}$};


%
%
%


\begin{scope}[shift={(60mm,-6.5mm)},scale=1]
\draw (0mm,0mm)\hezi;
\draw (5mm,0mm)\judy;
\draw (10mm,0mm)\hezi;
\draw (0mm,-5mm)\judy;
\draw (5mm,-5mm)\hezi;
\draw (10mm,-5mm)\hezi;
\draw (0mm,-10mm)\hezi;
\draw (5mm,-10mm)\hezi;
\draw (10mm,-10mm)\judy;


\node[rdot] at(5mm,0mm){};

\node at (7.5mm,-16mm){$213\in\mathcal{I}_{3,1}$};

\end{scope}
\end{tikzpicture}
\end{center}
Therefore, by the definitions of $\mathcal{B}^{(4)_3}_{n,k}$ and $\mathcal{B}^{(5)_3}_{n,k}$ given by \eqref{defB43nk} and \eqref{defB53nk},
we have
\begin{equation*}
  \left(\eta^{-1}_{i}(\pi),\left(i,\pi_{i}\right)\right)\in \mathcal{B}^{(4)_3}_{n-2,k}\uplus \mathcal{B}^{(5)_3}_{n-2,k-1}
  \subseteq \widetilde{\mathcal{B}}_{n,k},
\end{equation*}
since $\mathcal{B}^{(4)_3}_{n-2,k}\subseteq \mathcal{B}^{(3)}_{n-2,k}$ by \eqref{defB3nk},
and $\mathcal{B}^{(5)_3}_{n-2,k-1}\subseteq \mathcal{B}^{(4)}_{n-2,k-1}$ by \eqref{defB4nk}.

Consequently, for any $(\pi, i)\in\widetilde{\mathcal{I}}_{n,k}$,
we conclude that
\begin{equation*}
\Theta((\pi,i))\in\left\{
        \begin{array}{ll}
          \mathcal{B}^{(3)_1}_{n-2,k}\uplus\mathcal{B}^{(4)_1}_{n-2,k-1}
  \uplus\mathcal{B}^{(5)_1}_{n-2,k-2}, & \pi_i>i+1; \\[8pt]
          \mathcal{B}^{(3)_2}_{n-2,k}\uplus\mathcal{B}^{(4)_2}_{n-2,k-1}
  \uplus\mathcal{B}^{(5)_2}_{n-2,k-2}, & \pi_i<i-1; \\[8pt]
          \mathcal{B}^{(4)_3}_{n-2,k-1}\uplus\mathcal{B}^{(5)_3}_{n-2,k-2}, & \pi_i=i+1; \\[8pt]
          \mathcal{B}^{(2)_1}_{n-1,k-1}, & \pi_i=i-1; \\[8pt]
          \mathcal{B}^{(1)}_{n-1,k}\uplus \mathcal{B}^{(2)_2}_{n-1,k-1}, & \pi_i=i\text{ and $\pi_{i+1}\neq i+1$}; \\[8pt]
          \mathcal{B}^{(4)_3}_{n-2,k}\uplus \mathcal{B}^{(5)_3}_{n-2,k-1}, & \pi_i=i\text{ and $\pi_{i+1}=i+1$}.
        \end{array}
      \right.
\end{equation*}
On the other hand, we have
\begin{equation*}
\begin{aligned}
&\mathcal{B}^{(2)}_{n-1,k-1}=\mathcal{B}^{(2)_1}_{n-1,k-1}
\uplus\mathcal{B}^{(2)_2}_{n-1,k-1},\\[5pt]
&\mathcal{B}^{(3)}_{n-2,k}=\mathcal{B}^{(3)_1}_{n-2,k}
\uplus\mathcal{B}^{(3)_2}_{n-2,k}\uplus\mathcal{B}^{(4)_3}_{n-2,k},\\[5pt]
&\mathcal{B}^{(4)}_{n-2,k-1}=\mathcal{B}^{(4)_1}_{n-2,k-1}
\uplus\mathcal{B}^{(4)_2}_{n-2,k-1}\uplus\mathcal{B}^{(4)_3}_{n-2,k-1}
\uplus\mathcal{B}^{(5)_3}_{n-2,k-1},\\[5pt]
&\mathcal{B}^{(5)}_{n-2,k-2}=\mathcal{B}^{(5)_1}_{n-2,k-2}
\uplus\mathcal{B}^{(5)_2}_{n-2,k-2}\uplus\mathcal{B}^{(5)_3}_{n-2,k-2}.
\end{aligned}
\end{equation*}
Therefore, by the definition of $\widetilde{\mathcal{B}}_{n,k}$ given by \eqref{defBnk},
we derive that
\[
\Theta\colon \widetilde{\mathcal{I}}_{n,k}\rightarrow\widetilde{\mathcal{B}}_{n,k}
\]
is the inverse of $\Psi$,
which completes the proof.
\qed

\begin{exam}
For $\pi=42315\in\mathcal{I}_{5,2}$, we have
\begin{align*}
&\Theta(\pi,1)=(123,(1,3))\in\mathcal{B}^{(5)_1}_{3,0},
\Theta(\pi,2)=(213,(2,2))\in\mathcal{B}^{(5)_3}_{3,1},
\Theta(\pi,3)=(3214,(2,2))\in\mathcal{B}^{(1)}_{4,2},\\
&\Theta(\pi,4)=(123,(3,1))\in\mathcal{B}^{(5)_2}_{3,0},
\Theta(\pi,5)=(4231,(5,5))\in\mathcal{B}^{(2)_2}_{4,1},
\end{align*}
as listed below.
\begin{center}
\begin{tikzpicture}[scale = 0.65]
\def\hezi{-- +(5mm,0mm) -- +(5mm,5mm) -- +(0mm,5mm) -- cycle [line width=0.6pt]}
\def\judy{-- +(5mm,0mm) -- +(5mm,5mm) -- +(0mm,5mm) -- cycle [line width=0.6pt,fill=gainsboro]}
\tikzstyle{rdot}=[circle,fill=red,draw=red,inner sep=1.8]
\tikzstyle{cc}=[circle,draw=black,fill=yellow, line width=0.5pt, inner sep=1.5]

\draw[line width=0.8pt] (-5mm,33mm)--(185mm,33mm);
\node at (12.5mm, 40mm) {$\pi\in\mathcal{I}_{5,2}$};
\node at (50mm, 40mm) {$i=1$};
\node at (80mm, 40mm) {$i=2$};
\node at (110mm, 40mm) {$i=3$};
\node at (140mm, 40mm) {$i=4$};
\node at (170mm, 40mm) {$i=5$};
\draw[line width=1pt] (-5mm,47mm)--(185mm,47mm);


\draw[line width=0.8pt] (0mm,0mm)--(0mm,25mm);
\draw[line width=0.8pt] (5mm,0mm)--(5mm,25mm);
\draw[line width=0.8pt] (10mm,0mm)--(10mm,25mm);
\draw[line width=0.8pt] (15mm,0mm)--(15mm,25mm);
\draw[line width=0.8pt] (20mm,0mm)--(20mm,25mm);
\draw[line width=0.8pt] (25mm,0mm)--(25mm,25mm);

\draw[line width=0.8pt] (0mm,0mm)--(25mm,0mm);
\draw[line width=0.8pt] (0mm,5mm)--(25mm,5mm);
\draw[line width=0.8pt] (0mm,10mm)--(25mm,10mm);
\draw[line width=0.8pt] (0mm,15mm)--(25mm,15mm);
\draw[line width=0.8pt] (0mm,20mm)--(25mm,20mm);
\draw[line width=0.8pt] (0mm,25mm)--(25mm,25mm);

\draw (15mm,20mm)\judy;\draw (5mm,15mm)\judy;\draw (10mm,10mm)\judy;
\draw (0mm,5mm)\judy;\draw (20mm,0mm)\judy;
\node at (12.5mm, -5mm) {$42315$};

\begin{scope}[shift={(42.5mm,0mm)}]
\draw[line width=0.8pt] (0mm,0mm)--(0mm,15mm);
\draw[line width=0.8pt] (5mm,0mm)--(5mm,15mm);
\draw[line width=0.8pt] (10mm,0mm)--(10mm,15mm);
\draw[line width=0.8pt] (15mm,0mm)--(15mm,15mm);

\draw[line width=0.8pt] (0mm,0mm)--(15mm,0mm);
\draw[line width=0.8pt] (0mm,5mm)--(15mm,5mm);
\draw[line width=0.8pt] (0mm,10mm)--(15mm,10mm);
\draw[line width=0.8pt] (0mm,15mm)--(15mm,15mm);

\draw (0mm,10mm)\judy;
\draw (5mm,5mm)\judy;\draw (10mm,0mm)\judy;
\node[rdot] at (10mm,15mm){};

\node at (7.5mm, -5mm) {$123$};
\end{scope}

\begin{scope}[shift={(72.5mm,0mm)}]
\draw[line width=0.8pt] (0mm,0mm)--(0mm,15mm);
\draw[line width=0.8pt] (5mm,0mm)--(5mm,15mm);
\draw[line width=0.8pt] (10mm,0mm)--(10mm,15mm);
\draw[line width=0.8pt] (15mm,0mm)--(15mm,15mm);

\draw[line width=0.8pt] (0mm,0mm)--(15mm,0mm);
\draw[line width=0.8pt] (0mm,5mm)--(15mm,5mm);
\draw[line width=0.8pt] (0mm,10mm)--(15mm,10mm);
\draw[line width=0.8pt] (0mm,15mm)--(15mm,15mm);

\draw (5mm,10mm)\judy;
\draw (0mm,5mm)\judy;\draw (10mm,0mm)\judy;
\node[rdot] at (5mm,10mm){};

\node at (7.5mm, -5mm) {$213$};

\end{scope}

\begin{scope}[shift={(100mm,0mm)}]

\draw[line width=0.8pt] (0mm,0mm)--(0mm,20mm);
\draw[line width=0.8pt] (5mm,0mm)--(5mm,20mm);
\draw[line width=0.8pt] (10mm,0mm)--(10mm,20mm);
\draw[line width=0.8pt] (15mm,0mm)--(15mm,20mm);
\draw[line width=0.8pt] (20mm,0mm)--(20mm,20mm);

\draw[line width=0.8pt] (0mm,0mm)--(20mm,0mm);
\draw[line width=0.8pt] (0mm,5mm)--(20mm,5mm);
\draw[line width=0.8pt] (0mm,10mm)--(20mm,10mm);
\draw[line width=0.8pt] (0mm,15mm)--(20mm,15mm);
\draw[line width=0.8pt] (0mm,20mm)--(20mm,20mm);

\draw (10mm,15mm)\judy;\draw (5mm,10mm)\judy;
\draw (0mm,5mm)\judy;\draw (15mm,0mm)\judy;
\node[rdot] at (10mm,10mm){};
\node at (10mm, -5mm) {$3214$};

\end{scope}

\begin{scope}[shift={(132.5mm,0mm)}]
\draw[line width=0.8pt] (0mm,0mm)--(0mm,15mm);
\draw[line width=0.8pt] (5mm,0mm)--(5mm,15mm);
\draw[line width=0.8pt] (10mm,0mm)--(10mm,15mm);
\draw[line width=0.8pt] (15mm,0mm)--(15mm,15mm);

\draw[line width=0.8pt] (0mm,0mm)--(15mm,0mm);
\draw[line width=0.8pt] (0mm,5mm)--(15mm,5mm);
\draw[line width=0.8pt] (0mm,10mm)--(15mm,10mm);
\draw[line width=0.8pt] (0mm,15mm)--(15mm,15mm);

\draw (0mm,10mm)\judy;
\draw (5mm,5mm)\judy;\draw (10mm,0mm)\judy;
\node[rdot] at (0mm,5mm){};

\node at (7.5mm, -5mm) {$123$};
\end{scope}

\begin{scope}[shift={(160mm,0mm)}]
\draw[line width=0.8pt] (0mm,0mm)--(0mm,20mm);
\draw[line width=0.8pt] (5mm,0mm)--(5mm,20mm);
\draw[line width=0.8pt] (10mm,0mm)--(10mm,20mm);
\draw[line width=0.8pt] (15mm,0mm)--(15mm,20mm);
\draw[line width=0.8pt] (20mm,0mm)--(20mm,20mm);

\draw[line width=0.8pt] (0mm,0mm)--(20mm,0mm);
\draw[line width=0.8pt] (0mm,5mm)--(20mm,5mm);
\draw[line width=0.8pt] (0mm,10mm)--(20mm,10mm);
\draw[line width=0.8pt] (0mm,15mm)--(20mm,15mm);
\draw[line width=0.8pt] (0mm,20mm)--(20mm,20mm);

\draw (15mm,15mm)\judy;\draw (5mm,10mm)\judy;
\draw (10mm,5mm)\judy;\draw (0mm,0mm)\judy;
\node[rdot] at (20mm,0mm){};
\node at (10mm, -5mm) {$4231$};
\end{scope}

\draw[line width=0.6pt] (-5mm,-12mm)--(185mm,-12mm);

\draw[line width=0.4pt] (32mm,47mm)--(32mm,-12mm);
\draw[line width=0.4pt] (33mm,47mm)--(33mm,-12mm);
\end{tikzpicture}
\end{center}
\end{exam}

Let $\mathcal{J}_{n,k}$ be the set of all fixed-point free involutions on $[n]$ with $k$ descents.
By the restriction of fixed-point free on cycle decompositions \cite[Section 1.3]{Stanley-2012} of involutions,
we see $\mathcal{J}_{2n+1,k}=\emptyset$ for any $n,k\geq 0$.
Hence for fixed-point free involutions,
only the generating process
from the length of $(2n-2)$ to the length of $2n$ should be considered.

{\noindent{\emph{\textbf{Combinatorial Proof of Theorem \ref{thmJnk}.}}\hskip 2pt}}
Let
\begin{equation*}
  \widetilde{\mathcal{J}}_{n,k}=
  \{(\pi,i)\mid \pi\in \mathcal{J}_{n,k}\text{ and }1\leq i\leq n \}.
\end{equation*}
It is clear that $\left|\widetilde{\mathcal{J}}_{2n,k}\right|=2nJ_{2n,k}$,
which is the left side of \eqref{recJnk}.
We proceed to  construct several pairwise disjoint sets that
are counted by the terms on the right side of \eqref{recJnk}.

Let
\begin{equation}\label{defD1nkJ}
  \mathcal{D}^{(1)}_{n,k}=
  \{
  (\sigma,(i,j))\mid \sigma\in\mathcal{J}_{n,k}\text{ and $(i,j)$ is of $d$-type $(0,0)$ in $P_\sigma$ with $i\neq j$}
  \}.
\end{equation}
For $\sigma\in\mathcal{J}_{2n-2,k}$,
there are $((k+1)^2+2n-2)$ grid points of $d$-type $(0,0)$ in $P_\sigma$
by Theorem \ref{thmdtype}.
However, $(k+1)$ of such grid points lie on the main diagonal
since each $0_h$-path intersects with its symmetric $0_v$-path only once at the main diagonal due to the property of fixed-point free,
and there are $(k+1)$ $0_h$-paths by Theorem \ref{thmpaths}.
Thus we have
$$ \left|\mathcal{D}^{(1)}_{2n-2,k}\right|=(k(k+1)+2n-2)J_{2n-2,k},$$
which is the \textbf{first} term on the right side of \eqref{recJnk}.
For  $(\sigma, (i,j))\in \mathcal{D}^{(1)}_{2n-2,k}$, set
\begin{equation}\label{eqpfJn1st}
  \Psi((\sigma,(i,j)))=\left(\xi_{(i,j)}(\sigma),\chi_{ij}\right).
\end{equation}
Notice that the operations $\xi_{(i,j)}$ would not produce any new filled squares on the main diagonal,
which implies
$\left(\xi_{(i,j)}(\sigma),\chi_{ij}\right)\in\widetilde{\mathcal{J}}_{2n,k}$
by Proposition \ref{lemxi}.

Let
\begin{equation}\label{defD2nkJ}
  \mathcal{D}^{(2)}_{n,k}=
  \{
  (\sigma,(i,j))\mid \sigma\in\mathcal{J}_{n,k}\text{ and $(i,j)$ is of $d$-type $(1,0)$ or $(0,1)$ in $P_\sigma$}\},
\end{equation}
For $(\sigma,(i,j))\in \mathcal{D}^{(2)}_{2n-2,k-1}$,
set
\begin{equation}\label{eqpfJn2nd}
  \Psi((\sigma, (i,j)))=\left(\xi_{(i,j)}(\sigma),\chi_{ij}\right).
\end{equation}
By Proposition \ref{lemij=ji},
for all grid points in $\mathcal{D}^{(2)}_{2n-2,k-1}$, we have $i\neq j$.
And by Proposition \ref{lemxi},
we have $\left(\xi_{(i,j)}(\sigma),\chi_{ij}\right)\in\widetilde{\mathcal{J}}_{2n,k}$.

Let
\begin{equation}\label{defD3nkJ}
  \mathcal{D}^{(3)}_{n,k}=
  \{
  (\sigma,(i,i))_1,(\sigma,(i,i))_2\mid \sigma\in\mathcal{J}_{n,k} \text{ and $(i,i)$ is of $d$-type $(0,0)$ in $P_\sigma$}  \},
\end{equation}
where one can think that each pair $(\sigma,(i,i))$ in the set $\mathcal{D}^{(3)}_{n,k}$ appears twice,
and the subscripts $1$ and $2$ are used to indicate the first and second appearance, respectively.
For $r\in\{1,2\}$ and pairs $(\sigma, (i,i))_r\in \mathcal{D}^{(3)}_{2n-2,k-1}$,
let
\begin{equation}\label{eqpfJn3rd}
\Psi((\sigma,(i,i))_r)=\left\{
        \begin{array}{ll}
          \left(\eta'_i(\sigma),i\right), & r=1; \\[6pt]
          \left(\eta'_i(\sigma),i+1\right), & r=2.
        \end{array}
      \right.
\end{equation}
By 
Proposition \ref{lemetaeta'},
we obtain
$\left(\eta'_i(\sigma),i\right),\left(\eta'_i(\sigma),i+1\right)\in \widetilde{\mathcal{J}}_{2n,k}$.

Since there are $k$ grid points $(i,i)$ of $d$-type $(0,0)$ in $P_\sigma$ for $\sigma\in \mathcal{J}_{2n-2,k-1}$ by Theorem \ref{thmpaths},
together with Theorem \ref{thmdtype},
we get
\begin{align*}
  \left|\mathcal{D}^{(2)}_{2n-2,k-1}\uplus\mathcal{D}^{(3)}_{2n-2,k-1}\right|
  &=2(k(2n-1-k)-(2n-2))J_{2n-2,k-1}+2kJ_{2n-2,k-1} \\
  &=2((k-1)(2 n-k-1)+1) J_{2n-2,k-1},\end{align*}
which is the \textbf{second} term on the right side of \eqref{recJnk}.

%
%

Let
\begin{equation}\label{defD4nkJ}
  \mathcal{D}^{(4)}_{n,k}=
  \{
  (\sigma,(i,j))\mid \sigma\in\mathcal{J}_{n,k}\text{ and $(i,j)$ is of $d$-type $(1,1)$ in $P_\sigma$ with $i\neq j$}
  \},
\end{equation}
and
\begin{equation}\label{defD5nkJ}
  \mathcal{D}^{(5)}_{n,k}=
  \{
  (\sigma,(i,i))_1,(\sigma,(i,i))_2\mid \sigma\in\mathcal{J}_{n,k}\text{ and $(i,i)$ is of $d$-type $(1,1)$ in $P_\sigma$ }
  \}.
\end{equation}
For $(\sigma, (i,j))\in\mathcal{D}^{(4)}_{2n-2,k-2}$,
set
\begin{equation}\label{eqpfJn4th}
  \Psi((\sigma,(i,j)))=\left(\xi_{(i,j)}(\sigma),\chi_{ij}\right),
\end{equation}
and for $(\sigma, (i,i))_r\in\mathcal{D}^{(5)}_{2n-2,k-2}$ with $r\in\{1,2\}$, set
\begin{equation}\label{eqpfJn5th}
\Psi((\sigma,(i,i))_r)=\left\{
        \begin{array}{ll}
          \left(\eta'_i(\sigma),i\right), & r=1; \\[6pt]
          \left(\eta'_i(\sigma),i+1\right), & r=2.
        \end{array}
      \right.
\end{equation}
Following the similar analysis for \eqref{eqpfJn2nd} and \eqref{eqpfJn3rd},
we have $\left(\xi_{(i,j)}(\sigma),\chi_{ij}\right),\left(\eta'_i(\sigma),i\right),
\left(\eta'_i(\sigma),i+1\right)\in\widetilde{\mathcal{J}}_{2n,k}$.

For $\sigma\in \mathcal{J}_{2n-2,k-2}$,
it follows from Theorem \ref{thmdtype}
that there exist $((2n-k)^2+(2n-2))$ grid points with $d$-type $(1,1)$
in  $P_\sigma$.
Note that no squares on the main diagonal of $P_\sigma$ are filled,
which implies that each $1_h$-path should intersect with its symmetric
$1_v$-path exactly at the grid point on the main diagonal.
By theorem \ref{thmpaths},
we have $(2n-k)$ $1_h$-paths in  $P_\sigma$,
hence there are $(2n-k)$ grid points of $d$-type $(1,1)$ on the main diagonal.
Hence we see
\begin{align*}
\left|\mathcal{D}^{(4)}_{2n-2,k-2}\uplus\mathcal{D}^{(5)}_{2n-2,k-2}\right|
&=((2n-k)^2+2n-2-(2n-k))J_{2n-2,k-2}+2(2n-k)J_{2n-2,k-2}\\
&=((2 n-k)(2 n-k+1)+2 n-2) J_{2n-2,k-2},
\end{align*}
which is the \textbf{third} term on the right side of \eqref{recJnk}.
%
%

Therefore,  we construct five pairwise disjoint sets $\mathcal{D}^{(1)}_{2n-2,k}$,
$\mathcal{D}^{(2)}_{2n-2,k-1}$, $\mathcal{D}^{(3)}_{2n-2,k-1}$,
$\mathcal{D}^{(4)}_{2n-2,k-2}$, $\mathcal{D}^{(5)}_{2n-2,k-2}$
that together form the set
\begin{equation}\label{defD2nk}
  \widetilde{\mathcal{D}}_{2n,k}:=
  \mathcal{D}^{(1)}_{2n-2,k}\uplus\mathcal{D}^{(2)}_{2n-2,k-1}\uplus
  \mathcal{D}^{(3)}_{2n-2,k-1}\uplus\mathcal{D}^{(4)}_{2n-2,k-2}
  \uplus\mathcal{D}^{(5)}_{2n-2,k-2}
\end{equation}
satisfying that 
$\left|\widetilde{\mathcal{D}}_{2n,k}\right|$ equals  the right side of \eqref{recJnk}.
Furthermore, in terms of \eqref{eqpfJn1st}, \eqref{eqpfJn2nd}, \eqref{eqpfJn3rd}, \eqref{eqpfJn4th} and \eqref{eqpfJn5th},
we establish the mapping
\begin{equation*}
  \Psi\colon \widetilde{\mathcal{D}}_{2n,k}\rightarrow \widetilde{\mathcal{J}}_{2n,k}.
\end{equation*}

To complete the proof, for $(\pi,i)\in \widetilde{\mathcal{J}}_{2n,k}$,
we define
\begin{equation}\label{Jtheta}
\Theta((\pi,i))=\left\{
        \begin{array}{ll}
          \left(\xi^{-1}_{\langle i,\pi_i-1\rangle}(\pi),\left(i,\pi_i-1\right)\right), & \pi_i>i+1; \\[8pt]
          \left(\xi^{-1}_{\langle i-1,\pi_i\rangle}(\pi),\left(i-1,\pi_i\right)\right), & \pi_i<i-1; \\[8pt]
          \left(\eta'^{-1}_i(\pi),\left(i,\pi_i-1\right)\right)_1, & \pi_i=i+1; \\[8pt]
          \left(\eta'^{-1}_{i-1}(\pi),\left(i-1,\pi_i\right)\right)_2, & \pi_i=i-1.
        \end{array}
      \right.
\end{equation}

The cases of  $\pi_i>i+1$ and $\pi_i<i-1$ are exactly the same as
Cases $1$ and Cases  $2$ in the proof of Theorem \ref{thmInk}, respectively.
Let
\begin{align*}
  &\mathcal{D}^{(1)_1}_{n,k}=\{(\sigma,(i,j))\mid \sigma\in \mathcal{J}_{n,k}\text{, $(i,j)$ is of $d$-type $(0,0)$ in $P_\sigma$ with $i<j$}\},\\[3pt]
  &\mathcal{D}^{(2)_1}_{n,k}=\{(\sigma,(i,j))\mid \sigma\in \mathcal{J}_{n,k}\text{, $(i,j)$ is of $d$-type $(0,1)$ or $(1,0)$ in $P_\sigma$ with $i<j$}\}, \\[3pt]
  &\mathcal{D}^{(4)_1}_{n,k}=\{(\sigma,(i,j))\mid \sigma\in \mathcal{J}_{n,k}\text{, $(i,j)$ is of $d$-type $(1,1)$ in $P_\sigma$ with $i<j$}\},
\end{align*}
and
\begin{align*}
  &\mathcal{D}^{(1)_2}_{n,k}=\{(\sigma,(i,j))\mid \sigma\in \mathcal{J}_{n,k}\text{, $(i,j)$ is of $d$-type $(0,0)$ in $P_\sigma$ with $i>j$}\},\\[3pt]
  &\mathcal{D}^{(2)_2}_{n,k}=\{(\sigma,(i,j))\mid \sigma\in \mathcal{J}_{n,k}\text{, $(i,j)$ is of $d$-type $(0,1)$ or $(1,0)$ in $P_\sigma$ with $i>j$}\}, \\[3pt]
  &\mathcal{D}^{(4)_2}_{n,k}=\{(\sigma,(i,j))\mid \sigma\in \mathcal{J}_{n,k}\text{, $(i,j)$ is of $d$-type $(1,1)$ in $P_\sigma$ with $i>j$}\}.
\end{align*}
For $r\in\{1,2,4\}$ and $t\in\{1,2\}$, we deduce that $\mathcal{D}^{(r)_t}_{n,k}\subseteq \mathcal{D}^{(r)}_{n,k}$ in terms of \eqref{defD1nkJ}, \eqref{defD2nkJ} and \eqref{defD4nkJ}.
Thus for $\pi_i>i+1$, we have
\begin{equation*}
  \left(\xi^{-1}_{\langle i,\pi_i-1\rangle}(\pi),\left(i,\pi_i-1\right)\right)\in \mathcal{D}^{(1)_1}_{2n-2,k}\uplus\mathcal{D}^{(2)_1}_{2n-2,k-1}
  \uplus\mathcal{D}^{(4)_1}_{2n-2,k-2}
  \subseteq \widetilde{\mathcal{D}}_{2n,k},
\end{equation*}
and for $\pi_i<i-1$, we have
\begin{equation*}
  \left(\xi^{-1}_{\langle i-1,\pi_i\rangle}(\pi),\left(i-1,\pi_i\right)\right)\in \mathcal{D}^{(1)_2}_{2n-2,k}\uplus\mathcal{D}^{(2)_2}_{2n-2,k-1}
  \uplus\mathcal{D}^{(4)_2}_{2n-2,k-2}
  \subseteq \widetilde{\mathcal{D}}_{2n,k}.
\end{equation*}

If $\pi_i=i+1$, then both the squares $\langle i,i+1\rangle$ and $\langle i+1,i\rangle$ are filled in $P_\pi$.
Therefore, by Proposition \ref{lemetaeta'}, we have
$$\eta'^{-1}_i(\pi)\in \mathcal{J}_{2n-2,k-1-p}$$
and the grid point $(i,i)$ is of $d$-type $(p,p)$ in the grid of $\eta'^{-1}_i(\pi)$
for $p\in\{0,1\}$.
Hence, by setting
\begin{align*}
  &\mathcal{D}^{(3)_1}_{n,k}=\{(\sigma,(i,i))_1\mid \sigma\in \mathcal{J}_{n,k}\text{, $(i,i)$ is of $d$-type $(0,0)$ in $P_\sigma$}\},\\
  &\mathcal{D}^{(5)_1}_{n,k}=\{(\sigma,(i,i))_1\mid \sigma\in \mathcal{J}_{n,k}\text{, $(i,i)$ is of $d$-type $(1,1)$ in $P_\sigma$}\},
\end{align*}
we deduce that $\mathcal{D}^{(r)_1}_{n,k}\subseteq\mathcal{D}^{(r)}_{n,k}$ for $r\in\{3,5\}$ by \eqref{defD3nkJ} and \eqref{defD5nkJ},
and
\begin{equation*}
  \left(\eta'^{-1}_i(\pi),\left(i,\pi_i-1\right)\right)\in \mathcal{D}^{(3)_1}_{2n-2,k-1}\uplus\mathcal{D}^{(5)_1}_{2n-2,k-2}
  \subseteq \widetilde{\mathcal{D}}_{2n,k}.
\end{equation*}

If $\pi_i=i-1$, then the squares $\langle i-1, i\rangle$ and $\langle i,i-1\rangle $
are filled in $P_\pi$.
Thus based on the similar analysis of the case of $\pi_i=i-1$ above, by letting
\begin{align*}
  &\mathcal{D}^{(3)_2}_{n,k}=\{(\sigma,(i,i))_2\mid \sigma\in \mathcal{J}_{n,k}\text{, $(i,i)$ is of $d$-type $(0,0)$ in $P_\sigma$}\},\\
  &\mathcal{D}^{(5)_2}_{n,k}=\{(\sigma,(i,i))_2\mid \sigma\in \mathcal{J}_{n,k}\text{, $(i,i)$ is of $d$-type $(1,1)$ in $P_\sigma$}\},
\end{align*}
we have $\mathcal{D}^{(r)_2}_{n,k}\subseteq\mathcal{D}^{(r)}_{n,k}$ for $r\in\{3,5\}$ by \eqref{defD3nkJ} and \eqref{defD5nkJ},
and
\begin{equation*}
  \left(\eta'^{-1}_{i-1}(\pi),\left(i-1,\pi_i\right)\right)\in \mathcal{D}^{(3)_2}_{2n-2,k-1}\uplus\mathcal{D}^{(5)_2}_{2n-2,k-2}
  \subseteq \widetilde{\mathcal{D}}_{2n,k}.
\end{equation*}

In conclusion, for any pair $(\pi, i)$ in $\widetilde{\mathcal{J}}_{2n,k}$,
we show that
\begin{equation*}
\Theta((\pi,i))\in\left\{
        \begin{array}{ll}
         \mathcal{D}^{(1)_1}_{2n-2,k}\uplus\mathcal{D}^{(2)_1}_{2n-2,k-1}
  \uplus\mathcal{D}^{(4)_1}_{2n-2,k-2}, & \pi_i>i+1; \\[8pt]
       \mathcal{D}^{(1)_2}_{2n-2,k}\uplus\mathcal{D}^{(2)_2}_{2n-2,k-1}
  \uplus\mathcal{D}^{(4)_2}_{2n-2,k-2}, & \pi_i<i-1; \\[8pt]
       \mathcal{D}^{(3)_1}_{2n-2,k-1}\uplus\mathcal{D}^{(5)_1}_{2n-2,k-2}, & \pi_i=i+1; \\[8pt]
       \mathcal{D}^{(3)_2}_{2n-2,k-1}\uplus\mathcal{D}^{(5)_2}_{2n-2,k-2}, & \pi_i=i-1.
        \end{array}
      \right.
\end{equation*}
Particularly, we have
\begin{equation*}
\mathcal{D}^{(r)}_{n,k}=\mathcal{D}^{(r)_1}_{n,k}
\uplus\mathcal{D}^{(r)_2}_{n,k}
\end{equation*}
for $r\in \{1,2,3,4,5\}$.
Therefore, by \eqref{defD2nk}, we see that
\[
\Theta\colon \widetilde{\mathcal{J}}_{2n,k}\rightarrow\widetilde{\mathcal{D}}_{2n,k}
\]
is the inverse of $\Psi$.
\qed

\begin{exam}
For $\pi=532614\in\mathcal{J}_{6,3}$, we have
\begin{align*}
&\Theta(\pi,1)=(2143,(1,4))\in\mathcal{D}^{(2)_1}_{4,2},\quad
\Theta(\pi,2)=(3412,(2,2)_1)\in\mathcal{D}^{(5)_1}_{4,1},\\
&\Theta(\pi,3)=(3413,(2,2)_2)\in\mathcal{D}^{(5_2)}_{4,1},\quad
\Theta(\pi,4)=(4321,(4,5))\in\mathcal{D}^{(1)_1}_{4,3},\\
&\Theta(\pi,5)=(2143,(4,1))\in\mathcal{D}^{(2)_2}_{4,2},\quad
\Theta(\pi,6)=(4321,(5,4))\in\mathcal{D}^{(1)_2}_{4,3},
\end{align*}
as listed below.
\begin{center}
\begin{tikzpicture}[scale = 0.65]
\def\hezi{-- +(5mm,0mm) -- +(5mm,5mm) -- +(0mm,5mm) -- cycle [line width=0.6pt]}
\def\judy{-- +(5mm,0mm) -- +(5mm,5mm) -- +(0mm,5mm) -- cycle [line width=0.6pt,fill=gainsboro]}
\tikzstyle{rdot}=[circle,fill=red,draw=red,inner sep=1.8]
\tikzstyle{cc}=[circle,double=red,draw=red,fill=white, line width=0.5pt, inner sep=1.6]

\draw[line width=0.8pt] (-5mm,35mm)--(190mm,35mm);
\node at (15mm, 41mm) {$\pi\in\mathcal{J}_{6,3}$};
\node at (50mm, 41mm) {$i=1$};
\node at (75mm, 41mm) {$i=2$};
\node at (100mm, 41mm) {$i=3$};
\node at (125mm, 41mm) {$i=4$};
\node at (150mm, 41mm) {$i=5$};
\node at (175mm, 41mm) {$i=6$};
\draw[line width=1pt] (-5mm,47mm)--(190mm,47mm);


\draw[line width=0.8pt] (0mm,0mm)--(0mm,30mm);
\draw[line width=0.8pt] (5mm,0mm)--(5mm,30mm);
\draw[line width=0.8pt] (10mm,0mm)--(10mm,30mm);
\draw[line width=0.8pt] (15mm,0mm)--(15mm,30mm);
\draw[line width=0.8pt] (20mm,0mm)--(20mm,30mm);
\draw[line width=0.8pt] (25mm,0mm)--(25mm,30mm);
\draw[line width=0.8pt] (30mm,0mm)--(30mm,30mm);

\draw[line width=0.8pt] (0mm,0mm)--(30mm,0mm);
\draw[line width=0.8pt] (0mm,5mm)--(30mm,5mm);
\draw[line width=0.8pt] (0mm,10mm)--(30mm,10mm);
\draw[line width=0.8pt] (0mm,15mm)--(30mm,15mm);
\draw[line width=0.8pt] (0mm,20mm)--(30mm,20mm);
\draw[line width=0.8pt] (0mm,25mm)--(30mm,25mm);
\draw[line width=0.8pt] (0mm,30mm)--(30mm,30mm);

\draw (20mm,25mm)\judy;\draw (10mm,20mm)\judy;\draw (5mm,15mm)\judy;
\draw (25mm,10mm)\judy;\draw (0mm,5mm)\judy;\draw (15mm,0mm)\judy;

\node at (15mm, -5mm) {$532614$};

\begin{scope}[shift={(40mm,5mm)}]
\draw[line width=0.8pt] (0mm,0mm)--(0mm,20mm);
\draw[line width=0.8pt] (5mm,0mm)--(5mm,20mm);
\draw[line width=0.8pt] (10mm,0mm)--(10mm,20mm);
\draw[line width=0.8pt] (15mm,0mm)--(15mm,20mm);
\draw[line width=0.8pt] (20mm,0mm)--(20mm,20mm);

\draw[line width=0.8pt] (0mm,0mm)--(20mm,0mm);
\draw[line width=0.8pt] (0mm,5mm)--(20mm,5mm);
\draw[line width=0.8pt] (0mm,10mm)--(20mm,10mm);
\draw[line width=0.8pt] (0mm,15mm)--(20mm,15mm);
\draw[line width=0.8pt] (0mm,20mm)--(20mm,20mm);

\draw (5mm,15mm)\judy;\draw (0mm,10mm)\judy;
\draw (15mm,5mm)\judy;\draw (10mm,0mm)\judy;
\node[rdot] at (15mm,20mm){};

\node at (10mm, -5mm) {$2143$};

\end{scope}

\begin{scope}[shift={(65mm,5mm)}]

\draw[line width=0.8pt] (0mm,0mm)--(0mm,20mm);
\draw[line width=0.8pt] (5mm,0mm)--(5mm,20mm);
\draw[line width=0.8pt] (10mm,0mm)--(10mm,20mm);
\draw[line width=0.8pt] (15mm,0mm)--(15mm,20mm);
\draw[line width=0.8pt] (20mm,0mm)--(20mm,20mm);

\draw[line width=0.8pt] (0mm,0mm)--(20mm,0mm);
\draw[line width=0.8pt] (0mm,5mm)--(20mm,5mm);
\draw[line width=0.8pt] (0mm,10mm)--(20mm,10mm);
\draw[line width=0.8pt] (0mm,15mm)--(20mm,15mm);
\draw[line width=0.8pt] (0mm,20mm)--(20mm,20mm);

\draw (10mm,15mm)\judy;\draw (15mm,10mm)\judy;
\draw (0mm,5mm)\judy;\draw (5mm,0mm)\judy;
\node[rdot] at (5mm,15mm){};

\node at (10mm, -5mm) {$3412$};

\end{scope}

\begin{scope}[shift={(90mm,5mm)}]
\draw[line width=0.8pt] (0mm,0mm)--(0mm,20mm);
\draw[line width=0.8pt] (5mm,0mm)--(5mm,20mm);
\draw[line width=0.8pt] (10mm,0mm)--(10mm,20mm);
\draw[line width=0.8pt] (15mm,0mm)--(15mm,20mm);
\draw[line width=0.8pt] (20mm,0mm)--(20mm,20mm);

\draw[line width=0.8pt] (0mm,0mm)--(20mm,0mm);
\draw[line width=0.8pt] (0mm,5mm)--(20mm,5mm);
\draw[line width=0.8pt] (0mm,10mm)--(20mm,10mm);
\draw[line width=0.8pt] (0mm,15mm)--(20mm,15mm);
\draw[line width=0.8pt] (0mm,20mm)--(20mm,20mm);

\draw (10mm,15mm)\judy;\draw (15mm,10mm)\judy;
\draw (0mm,5mm)\judy;\draw (5mm,0mm)\judy;
\node[cc] at (5mm,15mm){};

\node at (10mm, -5mm) {$3412$};

\end{scope}

\begin{scope}[shift={(115mm,5mm)}]
\draw[line width=0.8pt] (0mm,0mm)--(0mm,20mm);
\draw[line width=0.8pt] (5mm,0mm)--(5mm,20mm);
\draw[line width=0.8pt] (10mm,0mm)--(10mm,20mm);
\draw[line width=0.8pt] (15mm,0mm)--(15mm,20mm);
\draw[line width=0.8pt] (20mm,0mm)--(20mm,20mm);

\draw[line width=0.8pt] (0mm,0mm)--(20mm,0mm);
\draw[line width=0.8pt] (0mm,5mm)--(20mm,5mm);
\draw[line width=0.8pt] (0mm,10mm)--(20mm,10mm);
\draw[line width=0.8pt] (0mm,15mm)--(20mm,15mm);
\draw[line width=0.8pt] (0mm,20mm)--(20mm,20mm);

\draw (15mm,15mm)\judy;\draw (10mm,10mm)\judy;
\draw (5mm,5mm)\judy;\draw (0mm,0mm)\judy;
\node[rdot] at (20mm,5mm){};

\node at (10mm, -5mm) {$4321$};
\end{scope}

\begin{scope}[shift={(140mm,5mm)}]
\draw[line width=0.8pt] (0mm,0mm)--(0mm,20mm);
\draw[line width=0.8pt] (5mm,0mm)--(5mm,20mm);
\draw[line width=0.8pt] (10mm,0mm)--(10mm,20mm);
\draw[line width=0.8pt] (15mm,0mm)--(15mm,20mm);
\draw[line width=0.8pt] (20mm,0mm)--(20mm,20mm);

\draw[line width=0.8pt] (0mm,0mm)--(20mm,0mm);
\draw[line width=0.8pt] (0mm,5mm)--(20mm,5mm);
\draw[line width=0.8pt] (0mm,10mm)--(20mm,10mm);
\draw[line width=0.8pt] (0mm,15mm)--(20mm,15mm);
\draw[line width=0.8pt] (0mm,20mm)--(20mm,20mm);

\draw (5mm,15mm)\judy;\draw (0mm,10mm)\judy;
\draw (15mm,5mm)\judy;\draw (10mm,0mm)\judy;
\node[rdot] at (0mm,5mm){};

\node at (10mm, -5mm) {$2143$};

\end{scope}

\begin{scope}[shift={(165mm,5mm)}]
\draw[line width=0.8pt] (0mm,0mm)--(0mm,20mm);
\draw[line width=0.8pt] (5mm,0mm)--(5mm,20mm);
\draw[line width=0.8pt] (10mm,0mm)--(10mm,20mm);
\draw[line width=0.8pt] (15mm,0mm)--(15mm,20mm);
\draw[line width=0.8pt] (20mm,0mm)--(20mm,20mm);

\draw[line width=0.8pt] (0mm,0mm)--(20mm,0mm);
\draw[line width=0.8pt] (0mm,5mm)--(20mm,5mm);
\draw[line width=0.8pt] (0mm,10mm)--(20mm,10mm);
\draw[line width=0.8pt] (0mm,15mm)--(20mm,15mm);
\draw[line width=0.8pt] (0mm,20mm)--(20mm,20mm);

\draw (15mm,15mm)\judy;\draw (10mm,10mm)\judy;
\draw (5mm,5mm)\judy;\draw (0mm,0mm)\judy;
\node[rdot] at (15mm,0mm){};

\node at (10mm, -5mm) {$4321$};

\end{scope}

\draw[line width=0.6pt] (-5mm,-10mm)--(190mm,-10mm);

\draw[line width=0.4pt] (35mm,47mm)--(35mm,-10mm);
\draw[line width=0.4pt] (36mm,47mm)--(36mm,-10mm);

\end{tikzpicture}
\end{center}
\end{exam}

At the end of this paper,
we remark that the statistics $\des$ in the symmetric group $\mathfrak{S}_n$ can be generalized to
the hyperoctahedral group $\mathfrak{B}_n$ as the statistics $\des^B$ and $\des_B$, see \cite{Moustakas-2019}.
In \cite{Gao-2022} and \cite{Gao-2024}, the authors verified that the geometric tools in this paper can be
generalized to the hyperoctahedral group,
and utilized to give combinatorial proofs for several recursive formulas
related to the joint distributions of descents and idescents on
signed permutations or singed involutions in $\mathfrak{B}_n$.

\vspace{2ex}

\noindent{\bf Acknowledgements}
This work is supported by the National Natural Science Foundation of China (12001078),
and the Natural Science Foundation of Chongqing (CSTB2022NSCQ-MSX0465).

%
%

\vspace{4ex}

{\footnotesize
\textsc {School of Science,  Chongqing University
of Posts and Telecommunications, Chongqing 400065, People’s Republic of China}

\emph{Email address:}
\tt{zkli@cqupt.edu.cn, xhliu7@163.com}
}

\end{document}